\documentstyle{article}
\begin{document}
\newtheorem{theorem}{Theorem}[section]
\newtheorem{corollary}[theorem]{Corollary}
\newtheorem{lemma}[theorem]{Lemma}
\newtheorem{proposition}[theorem]{Proposition}
\newtheorem{step}[theorem]{Step}
\newtheorem{example}[theorem]{Example}
\newtheorem{remark}[theorem]{Remark}

\font\sixbb=msbm6
\font\eightbb=msbm8
\font\twelvebb=msbm10 scaled 1095
\newfam\bbfam
\textfont\bbfam=\twelvebb \scriptfont\bbfam=\eightbb
                           \scriptscriptfont\bbfam=\sixbb

\def\bb{\fam\bbfam\twelvebb}
\newcommand{\Rea}{{\bb R}}
\newcommand{\Com}{{\bb C}}
\newcommand{\Int}{{\bb Z}}
\newcommand{\D}{{\bb D}}
\newcommand{\E}{{\bb E}}
\newcommand{\F}{{\bb F}}
\newcommand{\PP}{{\bb P}}
\newcommand{\card}{\rm card}
\newcommand{\calD}{\mathcal{D}}
\newcommand{\FF}{{\bb F}}
\newcommand{\rk}{{\rm rank}}
\newcommand{\Ker}{{\rm Ker}}
\renewcommand{\Im}{{\rm Im}}
\newcommand{\enp}{\begin{flushright} $\Box$ \end{flushright}}
\def\cD{{\mathcal{D}}}

\def\Sn{{\cal S}_n(F)}
\def\PSn{P{\cal S}_n(F)}
\def\cM{l}
\def\vp{\varphi}
\def\GL{\mbox{\rm GL}}
\def\LRa{\Leftrightarrow}
\def\bp{\noindent{\bf Proof. } }
\def\ep{\hfill $\Box$}

\def\N{\mathbb{{N}}}
\def\H{{\cal H}_n(D)}
\def\inv{\overline{\phantom{e}}}

\title{Hua's fundamental theorem of geometry of rectangular matrices over EAS division rings
\thanks{The second author was
supported by a grant from ARRS, Slovenia}}
\author{Cl\' ement de Seguins Pazzis\footnote{Laboratoire de Math\' ematiques de Versailles,
Universit\' e de Versailles Saint-Quentin-en-Yvelines,
45 avenue des Etats-Unis, 78035 Versailles cedex, France. dsp.prof@gmail.com}, \quad
Peter \v Semrl\footnote{Faculty of Mathematics and Physics, University of Ljubljana,
        Jadranska 19, SI-1000 Ljubljana, Slovenia. peter.semrl@fmf.uni-lj.si}
        }

\date{}
\maketitle

\begin{abstract}
The fundamental theorem of geometry of rectangular
matrices describes
the general form of bijective maps on the space of all $m\times n$
matrices over a division ring $\D$ which preserve adjacency in both directions.
This result proved by Hua in the nineteen forties has been recently improved in
several directions. One can study such maps without the bijectivity assumption or
one can try to get the same conclusion under the weaker assumption that
adjacency is preserved in one direction only. And the last possibility is to
study maps acting between matrix spaces of different sizes.
The optimal result would describe maps preserving adjacency in
one direction only acting between spaces of rectangular
matrices of different sizes in the absence of any regularity condition
(injectivity or surjectivity).

A division ring is said to be EAS if it is not isomorphic to any proper subring.
It has been known before that it is possible to construct adjacency preserving maps with
wild behaviour on matrices over division rings that are not EAS. For matrices over EAS
division rings it has been recently proved that adjacency preserving maps acting between
matrix spaces of different sizes satisfying a certain weak surjectivity condition
are either degenerate or of the expected simple standard form.
We will  remove this weak surjectivity assumption, thus solving completely the long
standing open problem of the optimal version of Hua's theorem.
\end{abstract}
\maketitle



\section{Introduction and statement of the main result}

Let $\D$ be a division ring and $m, n$ be positive integers. By $M_{m\times n} (\D )$ we denote
the set of all $m\times n$ matrices over $\D$.   When $m=n$ we write
$M_n (\D ) = M_{n\times n} (\D )$.
 Matrices $A,B \in M_{m\times n} (\D )$ are said to be adjacent
if ${\rm rank}\, (A-B) = 1$.

Let ${\cal V}$ be a
space of matrices. Recall that a  map $\phi : {\cal V} \to {\cal V}$
preserves adjacency in both directions if for every pair $A,B \in {\cal V}$ the
matrices $\phi (A)$ and $\phi (B)$ are adjacent if and only if $A$ and $B$ are
adjacent. We say that a map $\phi : {\cal V} \to {\cal V}$ preserves adjacency (in one
direction only) if $\phi (A)$ and $\phi(B)$ are adjacent whenever $A,B \in {\cal V}$
are adjacent. The study of such maps was initiated by Hua
in the series of papers \cite{Hu1} - \cite{Hu8}.

Hua's fundamental theorem of geometry of rectangular matrices
(see \cite{Wan}) states that for every bijective map $\phi : M_{m\times n}(\D) \to
M_{m\times n}(\D)$, $m,n \ge 2$, preserving adjacency in both directions there exist
invertible matrices $T\in M_{m}(\D)$, $S\in M_n (\D)$, a
matrix $R\in M_{m\times n}(\D)$, and an automorphism $\tau$ of the division ring $\D$
such that
\begin{equation}\label{hua1}
\phi ( A) = T A^\tau  S + R, \ \ \ A \in M_{m\times n} (\D).
\end{equation}
Here,
$A^\tau = [a_{ij}]^\tau = [\tau (a_{ij})]$ is a matrix obtained from $A$ by applying
$\tau$ entrywise.
In the square case $m=n$ we have the additional possibility
\begin{equation}\label{hua2}
\phi ( A) = T\  ^{t}(A^\sigma)  S + R, \ \ \ A \in M_n (\D),
\end{equation}
where $T,S,R$ are matrices in $M_n (\D)$ with $T,S$ invertible, $\sigma: \D \to \D$ is an
anti-automorphism, and $\,  ^{t}A$ denotes the transpose of $A$.
Clearly, the converse statement is true as well, that is, any map of
the form (\ref{hua1}) or (\ref{hua2}) is bijective and preserves adjacency in both directions.

Composing the map $\phi$ with a translation affects neither the assumptions, nor
the conclusion of Hua's theorem. Thus, there is no loss of generality in assuming that $\phi(0) = 0$. Then clearly, $R=0$.
It is a remarkable fact that
after this harmless normalization the additive (semilinear in the case when $\D$ is a field) character of $\phi$ is not an
assumption but a conclusion.

Besides Hua and Siegel's original motivation, related respectively
to complex analysis and symplectic geometry, this result has many other applications, say in the theory of linear preservers \cite{LiP} and
in the geometry of Grassmann spaces. We refer an interested reader to \cite{Se6}
for a detailed explanation. In particular, it is expected that our main result and proof techniques
will lead to an improvement of Chow's fundamental theorem of geometry of Grassmann
spaces \cite{Cho}.

In Hua's theorem we get a very strong conclusion under very mild assumptions.
Nevertheless, one can ask if it is optimal. Or is it possible to improve it?
The first natural question is whether one can still get the same conclusion in replacing
the assumption that adjacency is preserved in both directions by the weaker assumption that
it is preserved in one direction only. This question had been
open for quite a long time  and then finally answered in the affirmative in \cite{HuW}.
One can further ask whether it is possible to relax the bijectivity assumption.
The first guess that Hua's theorem remains valid without the bijectivity assumption
with a minor modification that $\tau$ appearing in (\ref{hua1}) is an endomorphism of $\D$
(not necessarily surjective), while $\sigma$ appearing in (\ref{hua2}) is an
anti-endomorphism, turned out to be wrong for matrices over general division rings.
In particular, this conjecture is true for real matrices
and wrong for complex matrices \cite{Se1}.
And finally, one can try to improve Hua's theorem by studying
adjacency preserving maps  acting between spaces of matrices of
different sizes.

It should be mentioned that the analogous problems for adjacency preserving maps
on $H_n$, the space of all $n\times n$ complex hermitian matrices, are much easier.
In this case the Hua's theorem states that
every bijective map $\phi : H_n \to H_n$ preserving adjacency in both directions
and satisfying $\phi (0) = 0$ is
a congruence transformation (that is, a map of the form $M \mapsto P^\star MP$ with $P$ an invertible matrix)
possibly composed with the transposition and possibly
multiplied by $-1$. And again we can ask for possible improvements in all three
above mentioned directions. All three questions have been answered
simultaneously in the paper \cite{HuS} by obtaining the following optimal
result. Let $m,n$ be integers with $m\ge 2$ and $\phi : H_m \to H_n$ be a map
preserving adjacency (in one direction only; note that no surjectivity or injectivity is assumed
and that $m$ may be different from $n$) and satisfying $\phi (0) =0$ (this is,
of course, a harmless normalization). Then either $\phi$ is the
standard embedding of $H_m$ into $H_n$ composed with a congruence
transformation on $H_n$ possibly composed with the transposition and possibly multiplied
by $-1$; or $\phi$ is of a very special degenerate form, that is, its range is contained
in the linear span of some rank one hermitian matrix. This result has already been proved to be useful
including some applications in mathematical physics \cite{Se3,Se4,Se9}.

In the case of full matrix spaces the following theorem was proved
in \cite{Se6}: if
$\D$ is any division ring, ${\rm card}\, \D > 3$, $n,p,q$ are integers with $p,q \ge n \ge 3$, a map
 $\phi : M_{n} (\D) \to M_{p \times q} (\D)$ preserves adjacency and satisfies $\phi (0) =0$, and there exists
$A_0 \in M_{n} (\D)$ such that ${\rm rank}\, \phi (A_0 ) = n$,
then either there exist invertible matrices $T\in M_p (\D)$ and $S\in M_q (\D)$,
an endomorphism $\tau : \D \to \D$, and a matrix $L\in M_n (\D)$ with the property
that $I + A^\tau L \in M_n (\D)$ is invertible for every $A\in M_n (\D)$, such that
\begin{equation}\label{puksicko}
\phi (A) = T \left[ \begin{matrix}{ (I + A^\tau L)^{-1} A^\tau & 0 \cr 0 & 0 \cr}\end{matrix}\right] S, \ \ \ A\in M_n (\D);
\end{equation}
or  there exist invertible matrices $T\in M_p (\D)$ and $S\in M_q (\D)$, an
anti-endomorphism $\sigma : \D \to \D$, and a matrix $L\in M_n (\D)$ with the property
that $I + \, ^{t}(A^\sigma) L \in M_n (\D)$ is invertible for every $A\in M_n (\D)$, such that
\begin{equation}\label{gremoprot}
\phi (A) = T \left[ \begin{matrix}{ (I + \, ^{t}(A^\sigma) L)^{-1} \, ^{t}(A^\sigma) & 0 \cr 0 & 0 \cr}\end{matrix}\right] S, \ \ \ A\in M_n (\D);
\end{equation}
or $\phi$ is of some very special degenerate form.

Note that in this theorem the domain of $\phi$ is the space of square $n\times n$ matrices. In \cite{Se6} examples were given
showing that there is no nice description of adjacency preserving maps acting on a space of rectangular matrices and that the
assumption on the existence of
$A_0 \in M_{n} (\D)$ satisfying ${\rm rank}\, \phi (A_0 ) = n$ is indispensable, while the assumption that $n \ge 3$ was
needed because of the proof techniques and we believe it is superfluous.
The construction of the examples
was based on the existence of a nonsurjective endomorphism of the underlying division ring $\D$.

Let us recall that a division ring $\D$ is said to be an EAS\footnote{EAS stands for: ``Endomorphisms Are Surjective".} division ring if every
endomorphism $\tau : \D \to \D$
is automatically surjective. The field of real numbers and the field of rational numbers are
well-known to be EAS. Obviously, every finite field is EAS.
The same is true for the division ring of quaternions  (see, for example \cite{Semr}),
 while the complex field  is not an EAS field \cite{Kes}.
Let $\D$ be an EAS division ring.
 It is then easy to verify that also each anti-endomorphism of $\D$
is bijective (just note that the square of an anti-endomorphism is an
endomorphism).

The above theorem has a very simple form under the additional assumption that $\D$ is
an EAS division ring. Indeed, in this case every endomorphism $\tau : \D \to \D$ is an automorphism, and therefore,
 $I + A^\tau L \in M_n (\D)$ is invertible for every $A\in M_n (\D)$ if and only if $L=0$, and thus the formula
(\ref{puksicko}) has the following very simple form in this special case:
$$
\phi (A) = T \left[ \begin{matrix}{ A^\tau & 0 \cr 0 & 0 \cr}\end{matrix}\right] S, \ \ \ A\in M_n (\D).
$$
The equation (\ref{gremoprot}) becomes similarly simple under the EAS assumption.

It is now time to recall the notion
of standard adjacency preserving maps.
Let $m,n,p,q$ be positive integers with $p\ge m$ and $q\ge n$,
$\tau: \D \to \D$ an automorphism, $T \in M_p (\D)$ and $S \in M_q (\D )$ invertible
matrices, and $R\in M_{p \times q}(\D)$ any matrix. Then the map $\phi : M_{m\times n}(\D)
\to M_{p \times q} (\D)$ defined by
\begin{equation}\label{onopr}
\phi (A ) = T \left[ \begin{matrix} { A^\tau & 0\cr 0&0\cr }  \end{matrix} \right] S + R
\end{equation}
preserves adjacency in both directions.
Similarly, if  $m,n,p,q$ are positive integers with $p\ge n$ and $q\ge m$,
$\sigma: \D \to \D$ an anti-automorphism, $T \in M_p (\D)$ and $S \in M_q (\D )$ invertible
matrices, and $R \in M_{p \times q}(\D)$ any matrix, then $\phi : M_{m\times n}(\D)
\to M_{p \times q} (\D)$ defined by
\begin{equation}\label{onodr}
\phi (A ) = T \left[ \begin{matrix} {  \  ^{t}(A^\sigma) & 0\cr 0&0\cr }  \end{matrix} \right] S + R
\end{equation}
preserves adjacency in both directions as well. We will call any map that is of one of the above two forms a standard
adjacency preserving map.

Using the above theorem on adjacency preserving maps on square matrices over general division rings the following
corollary was obtained in \cite{Se6}.
Let $m,n,p,q$ be integers with $m,p,q , n \ge 3$ and $\D$ an EAS division ring, ${\rm card}\, \D > 3$. Assume that
 $\phi : M_{m \times n} (\D) \to M_{p \times q} (\D)$ preserves adjacency, $\phi (0) =0$, and there exists
$A_0 \in M_{m \times n} (\D)$ such that ${\rm rank}\, \phi (A_0 ) = \min \{ m,n \}$.
Then either $\phi$ is of a standard form, or it is of some rather special degenerate form.

The aim of this paper is to obtain the optimal version of Hua's theorem for matrices over
EAS division rings. In order to do this we have to remove both the assumption on the existence
of a matrix $A_0$ with ${\rm rank}\, \phi (A_0 ) = \min \{ m,n \}$ and the assumption that
$m,n,p,q \ge 3$. We will need a completely different approach from the one in \cite{Se6}. There the main idea
was to reduce the problem to the study of order preserving maps on idempotent matrices. Here
we will start with the characterization of adjacency preserving maps on 2 by 2 matrices. To obtain the general case we then need
to put the information on small pieces together to get a global picture.

There has been only one notion left to be clarified before we formulate our main result.
When speaking of known results we have mentioned adjacency preserving maps of some rather special degenerate form
several times. Recall first that a set ${\cal S}$ of matrices is said to be adjacent if every pair of matrices
$A,B \in {\cal S}$, $A\not=B$, is adjacent. As an example of such a set we can take a subset $\{ tR \, : \, t \in \Rea \}
\subset H_n$, where $R$ is a rank one hermitian matrix. Indeed, ${\rm rank}\, (tR-sR) = 1$ whenever $t\not=s$.
Hence, the above mentioned optimal version of Hua's fundamental theorem of geometry of hermitian matrices states
that every adjacency preserving map from $H_m$ into $H_n$, $m\ge 2$,
is either of a standard form, or
its image is contained in an adjacent set. With this result in mind one can ask whether
the assumption that $\D$ is an EAS division ring guarantees that
each adjacency preserving map $\phi : M_{m \times n} (\D ) \to M_{p \times q} (\D)$
is either standard, or its image is contained in some adjacent set.
The answer is negative as the following example from $\cite{Se6}$ shows.

\begin{example}\label{ggzzgg}
Let
$m \ge 4$, $\D$ be an infinite division ring, and
$$
\varphi_k : M_{m \times 4}^k  (\D) \to \D, \ \ \ k=1,2,3,4,
$$
be injective maps. Here, $M_{m \times 4}^k  (\D)$ denotes the set of all $m \times 4$ matrices of rank $k$.
Assume further that the ranges of $\varphi_3$ and $\varphi_4$ are disjoint.
Denote by $E_{i,j} \in M_4 (\D)$ the matrix units, that is, $E_{i,j}$ is the matrix whose entries are all zero but the $(i,j)$-entry that
is equal to $1$. Define the map $\phi : M_{m \times 4} (\D) \to M_4 (\D)$ by
$$
\phi (0) = 0,
$$
$$
\phi (A) = E_{1,1} + \varphi_1 (A) E_{1,2}, \ \ \ A \in M_{m \times 4}^1 (\D),
$$
$$
\phi (A) = E_{1,1} + E_{2,2} +  \varphi_2 (A) E_{3,2}, \ \ \ A \in M_{m \times 4}^2 (\D),
$$
$$
\phi (A) = E_{1,1} + E_{2,2} + E_{3,3} + \varphi_3 (A) E_{3,4}, \ \ \ A \in M_{m \times 4}^3 (\D),
$$
and
$$
\phi (A) = E_{1,1} + E_{2,2} + E_{3,3} + \varphi_4 (A) E_{3,4}, \ \ \ A \in M_{m \times 4}^4 (\D).
$$
It is straightforward to check that this map preserves adjacency. Obviously, $\phi$ is neither standard, nor its range is an adjacent set.
\end{example}

In the next section we will give a precise definition of the rank of a matrix over a division ring and then we will define the distance between
two matrices with the following formula
$$
\delta(A, B) = {\rm rank}\, (A-B), \ \ \ A,B \in M_{m\times n} (\D ).
$$
Then we can define the unit ball with a central point $A$ as ${\cal B} (A,1) = \{ B \in   M_{m\times n} (\D )\, : \, \delta(A,B) \le 1 \}$,
that is,  ${\cal B} (A,1) = A + M_{ m \times n}^{\le 1} (\D)$, where $M_{ m \times n}^{\le 1} (\D)$ denotes the set of all $m \times n$ matrices of rank at most one.

We need some more notation before we formulate the precise definition of a degenerate map.
Let $^{t}x$ and $y$ be any nonzero $p\times 1$ matrix and $1\times q$ matrix, respectively.
We further denote by $L(\, ^{t}x ) \subset M_{ p \times q} (\D)$ the set of all $p \times q$ matrices of the form $^{t}x u$, where $u$ is any
$1 \times q$ matrix. Clearly, $L (\, ^{t}x)$ is an adjacent set, and so is, $E+ L(\, ^{t}x ) =
\{ E + \, ^{t}x u \, : \, u \in M_{1 \times q} (\D) \}$. Here, $E$ is any $p \times q$ matrix. Similarly,
we define $R(y) = \{ \, ^{t}vy \, : \, \, ^{t}v \in M_{p \times 1} (\D) \}$.
An adjacency preserving map $\phi : M_{m \times n} (\D) \to M_{p \times q}(\D)$
is called degenerate with respect to a matrix $A \in M_{m \times n} (\D)$
if there exists a nonzero $p \times 1$ matrix
$^{t}x$ and a nonzero $1 \times q$ matrix $y$ such that
\begin{equation}\label{ahmlaj}
\phi ({\cal B} (A,1)) = \phi ( A + M_{m\times n}^{\le 1} (\D) ) \subset (\phi (A) + R(y)) \cup (\phi (A) + L (\, ^{t}x )).
\end{equation}
An adjacency preserving map $\phi : M_{m \times n} (\D) \to M_{p \times q}(\D)$
is called degenerate if it is degenerate with respect to every single matrix of $M_{m \times n} (\D)$.
In other words, a map $\phi$ is degenerate if it shrinks every unit ball (with respect to the arithmetic distance $\delta$)
into a very small set - a union of two adjacent sets. One can easily verify that this condition is satisfied for the map $\phi$ defined in Example \ref{ggzzgg}.
Every map $\phi : M_{m \times n} (\D) \to M_{p \times q} (\D)$ satisfying (\ref{ahmlaj}) for all $A\in M_{m \times n} (\D)$ almost preserves adjacency. Namely, if $\phi$ has this property and $B,C \in M_{m\times n}(\D)$ are adjacent, then either $\phi (B) = \phi (C)$, or $\phi (B)$ and $\phi (C)$ are adjacent.
An interested reader can find several examples of degenerate adjacency preserving maps in \cite{Se6}.
We are now ready to formulate our main result.

\begin{theorem}\label{maintheorem}
Let $m,n,p,q$ be positive integers, $\D$ be an EAS division ring, and
$\phi : M_{m \times n} (\D) \to M_{p \times q}(\D)$ be an adjacency preserver.
Then either $\phi$ is degenerate or it is standard.
\end{theorem}

Having this theorem we can say that the problem of finding the optimal version of Hua's fundamental theorem
of geometry of rectangular matrices has been completely solved. The case of matrices over
general division rings has been studied in \cite{Se6}. The optimality of the main theorem
in \cite{Se6} has been illustrated by several examples most of them based on the existence
of nonsurjective endomorphisms of the underlying division ring. This left open
only the special case when $\D$ is EAS which has been now solved in an optimal way.

At this point, one may naturally wonder whether a more explicit description of degenerate
adjacency preservers can be obtained. For infinite division rings,
the many examples from \cite{Se6} suggest that this problem should be intractable,
although we believe that additional results on the local behavior of degenerate maps can be obtained.
In the case of finite fields, the following theorem can be regarded as an almost full answer to this question:

\begin{theorem}\label{finitefieldstheo}
Let $\F$ be a finite field, $m,n,p,q$ be positive integers, and
$\phi : M_{m \times n}(\F) \rightarrow M_{p \times q}(\F)$ be a degenerate adjacency preserver.
Then, the range of $\phi$ is an adjacent set.
\end{theorem}

The only remaining problem for a finite field $\F$ would be to find
the minimal cardinality among the sets $X$ such that there is a
map $\phi : M_{m \times n}(\F) \rightarrow X$ that takes adjacent matrices to distinct elements.
This seems to be a difficult question, and at this point we have no reasonable conjecture for its solution.

\section{Notation and basic definitions}

We start by recalling the definition of the rank of an $m\times n$ matrix $A$ with entries
in a division ring $\D$. We will always consider $\D^n$, the set of all $1\times n$
matrices,   as a left vector space over $\D$.
Correspondingly,  we have the right vector
space of all $m\times 1$ matrices $^{t}\D^m$.
The row space of $A$ is defined to be the left vector subspace of $\D^n$ generated
by the rows
of $A$,
and the row rank of $A$ is defined to be the dimension of this subspace.
Similarly, the column rank of $A$ is the dimension of the right vector space
generated by the columns of $A$.
This space is called the column space of $A$.
These two ranks are equal for every matrix over $\D$ and this common value is called
the rank of a matrix. If ${\rm rank}\, A = r$, then there exist invertible matrices
$T\in M_m (\D)$ and $S\in M_{n}(\D)$ such that
\begin{equation}\label{cxy}
TAS = \left[\matrix{I_r & 0 \cr 0 & 0 \cr}\right] .
\end{equation}
Here, $I_r$ denotes the $r\times r$ identity matrix and the zeroes stand for zero matrices of the
appropriate sizes.
For a positive integer $r$, $1 \le r \le \min \{m , n\}$, we denote by $M_{m \times n}^r  (\D)$
the set of all matrices $A\in M_{m \times n} (\D)$ of rank $r$. We write shortly $M_{n}^r  (\D) = M_{n \times n}^r  (\D)$.

Since rank satisfies the triangle inequality, that is, ${\rm rank}\, (A+B) \le {\rm rank}\, A
+ {\rm rank}\, B$ for every pair $A,B \in M_{m\times n} (\D )$ \cite[p.46, Exercise 2]{Jac},
the set of matrices $M_{m\times n} (\D)$ equipped with the distance $\delta$ defined by
$$
\delta(A, B) = {\rm rank}\, (A-B), \ \ \ A,B \in M_{m\times n} (\D ),
$$
is a metric space. The distance $\delta$ is called the arithmetic distance. Two matrices
$A$ and $B$ are adjacent if and only if $\delta(A,B) = 1$.
An important point that we need in the proof of our main theorem is that adjacency preserving maps are
contractions with respect to $\delta$.
Indeed, assume that  $\phi : M_{m\times n}(\D)
\to M_{p \times q} (\D)$ preserves adjacency.
Using the facts
that $\delta$ satisfies the triangle inequality and that for every positive integer $r$ and
every pair $A,B \in M_{m \times n} (\D)$ we have
$\delta(A,B) = r$ if and only if there exists a chain of matrices
$A = A_0 , A_1 , \ldots , A_r = B$ such that the pairs $A_0 , A_1$, and
$A_1 , A_2$, and $\ldots$, and $A_{r-1} , A_r$ are all adjacent we
easily see that
$$
\delta(\phi (A), \phi (B) ) \le d(A,B), \ \ \ A,B \in M_{m \times n} (\D).
$$

The rank of a matrix $A$ need not be equal to the rank of its transpose $ ^{t}A$.
Let $\tau : \D \to \D$ be an anti-endomorphism, that is, $\tau$ is additive, $\tau(1)=1$
and $\tau (\lambda \mu) = \tau (\mu) \tau (\lambda)$, $\lambda, \mu \in \D$,
and denote by  $A^\tau = [a_{ij}]^\tau = [\tau (a_{ij})]$ the matrix obtained from
$A$ by applying $\tau$ entry-wise.
Then
${\rm rank}\, A = {\rm rank} \  ^{t}(A^\tau )$.

Let $a \in \D^n$ and $^{t}b \in \  ^{t}\D^m$
be any nonzero vectors. Then $^{t}ba = ( ^{t}b)a$ is a matrix
of rank one. Every matrix of rank one can be written in this form.
It is easy to verify that two rank one matrices
$^{t}ba$ and $^{t}dc$, $^{t}ba \not=\,  ^{t}dc$, are adjacent if and only if
$a$ and $c$ are linearly dependent or $^tb$ and $^td$ are linearly dependent.

For a nonzero $x \in \D^n$ and a nonzero $^{t}y \in \, ^{t}\D^m$ we denote by $R(x)$
and $L (\, ^{t}y )$ the subsets of $M_{m\times n} (\D )$ defined by
$$
R (x) = \{ \, ^{t}ux \, : \ ^{t}u \in \, ^{t}\D^m  \}
$$
and
$$
L (\, ^{t}y) = \{ \, ^{t}yv \, : \, v \in \D^n  \} .
$$
A subset ${\cal S} \subset M_{m \times n} (\D )$ is called an adjacent set if any two
distinct matrices $A,B \in {\cal S}$ are adjacent. Clearly, the above two sets are adjacent.
It is well-known (and easy to verify) that they are maximal adjacent sets. And any adjacent set
containing $0$ is contained in one of such subsets. Let $A\in M_{m \times n} (\D)$. Obviously, a
subset ${\cal S} \subset M_{m \times n} (\D )$ is an adjacent set if and only if
$A + {\cal S} = \{ A + B \, : \, B \in {\cal S} \}$ is an adjacent set. Thus,
 ${\cal T} \subset M_{m \times n} (\D )$ is a maximal adjacent set if and
only if ${\cal T} = A + R(x)$ for some $A\in M_{m \times n} (\D)$
and some nonzero $x \in \D^n$, or
${\cal T} = A + L(\, ^{t}y)$ for some $A\in M_{m \times n} (\D)$
and some nonzero $\, ^{t}y \in \, ^{t}\D^m$.

Remember that $\D^n$ and ${}^t \D^m$ are respectively equipped with canonical structures of left and right vector spaces
over $\D$; in what follows we will denote by $\PP(\D^n)$ and $\PP({}^t \D^m)$ their respective projective spaces,
i.e.\ $\PP(\D^n)$ (respectively $\PP({}^t \D^m)$) is the set of all $1$-dimensional subspaces of $\D^n$
(respectively, of ${}^t\D^m$). Given a point $d \in \PP(\D^n)$ and a nonzero vector $x \in d$, one sees that
the adjacent set $R(x)$ is independent from the choice of $x$, and we shall denote it by $R_d$.
Similarly, given a point $e \in \PP({}^t\D^m)$, the set $L({}^t y)$ does not depend on the choice of a nonzero vector $y \in e$,
so that we may define $L_e:= L({}^t y)$.

As usual, the symbol $E_{i,j}$, $1\le i \le m$, $1 \le j \le n$, will stand for the matrix
having all entries zero but the $(i,j)$-entry which is equal to $1$. Later on we will deal simultaneously with
rectangular matrices of different sizes, say with matrices
from $M_{m \times n} (\D)$ and $M_{p \times q} (\D)$.
The same symbol $E_{i,j}$ will be used to denote a matrix
unit in $M_{m \times n} (\D)$ as well as a matrix
unit in $M_{p \times q} (\D)$.

The set of $m\times n$ matrices will be identified with the set of linear transformations mapping
$\D^m$ into $\D^n$. Each
$m \times n$ matrix $A$ gives rise to a linear operator
defined by $x\mapsto xA$, $x\in \D^m$. The rank of the matrix $A$ is equal to the dimension
of the image ${\rm Im}\, A$ of the corresponding operator $A$. The kernel of an operator $A$ is defined as
${\rm Ker}\, A = \{ x \in \D^m \, : \, xA = 0 \}$. We have $m = {\rm rank}\, A + \dim {\rm Ker}\, A$.

One of the main tools in the proof of our main result is a non-surjective version of the fundamental theorem of
affine geometry. Given vector spaces $U$ and $V$ over $\D$ (either left or right), a semilinear map from $U$ to $V$
is a group homomorphism $f$ from $U$ to $V$ for which there is a bijection $\sigma : \D \rightarrow \D$
which satisfies the set of conditions on one of the rows of the following array:

\begin{center}
\begin{tabular}{| c | c | c | c |}
\hline
$U$ & $V$  & $\sigma$ & For all $(\lambda,x) \in \D \times U$,  \\
\hline
\hline
left vector space & left vector space & automorphism& $f(\lambda\,x)=\lambda^\sigma f(x)$ \\
\hline
left vector space & right vector space & anti-automorphism & $f(\lambda\,x)=f(x) \lambda^\sigma$ \\
\hline
right vector space & left vector space & anti-automorphism & $f(x\,\lambda)=\lambda^\sigma\,f(x)$ \\
\hline
right vector space & right vector space & automorphism & $f(x\,\lambda)=f(x)\,\lambda^\sigma$ \\
\hline
\end{tabular}
\end{center}
A map from $U$ to $V$ is called a lineation if it maps every affine line into an affine line.

Assume that $\D$ is an EAS division ring with ${\rm card}\, \D > 2$.
Then, the fundamental theorem of affine geometry states that
every injective lineation from $U$ to $V$ sending $0$ to $0$ whose range is not included in an affine line is semilinear
(and even a semilinear isomorphism if the spaces $U$ and $V$ are finite-dimensional with the same dimension).
An interested reader can find the proof in the case when $\D$ is commutative in \cite[page 104]{Ben}.
The above version of the fundamental theorem of affine geometry is due to Schaeffer \cite{Sch} who
formulated and proved his result also for general (not necessarily commutative) division rings.

\section{Proof of the main result}

\subsection{Basic results on adjacent matrices}\label{basicresults}

Let $\D$ be a division ring, $m,n$ be positive integers, and $A, B \in M_{m\times n} (\D )$.
Assume that
\begin{equation}\label{mijegod}
{\rm rank}\, (A+B) = {\rm rank}\, A + {\rm rank}\, B .
\end{equation}
We identify matrices with operators. Then (see \cite{Se6})
$$
{\rm Im}\, (A+B) = {\rm Im}\, A \oplus {\rm Im}\, B
$$
and
$$
{\rm Ker}\, (A+B) = {\rm Ker}\, A \cap {\rm Ker}\, B.
$$
In particular, if $A$ and $B$ are adjacent and ${\rm rank}\, A < {\rm rank}\, B$, then $B= A + R$
for some $R$ of rank one. It follows that ${\rm rank}\, B = {\rm rank}\, (A+R) = {\rm rank}\, A +
{\rm rank}\, R$, and therefore, ${\rm Im}\, A \subset {\rm Im}\, B$ and ${\rm Ker}\, B \subset
{\rm Ker}\, A$.

The proofs of the following two lemmas can be found in \cite[Lemmas 5.1 and 5.2]{Se6}.

\begin{lemma}\label{subadjacentidempotent}
Let $n,p,q$ be positive integers with $n\le p,q$. Set
$$P=\left[ \matrix{ I_n & 0 \cr 0 & 0 \cr}\right] \in M_{p\times q} (\D).$$
Then, the rank $n-1$ matrices of $M_{p \times q}(\D)$ that are adjacent to $P$ are the matrices of the form
$$
Q = \left[ \matrix{ Q_1 & 0 \cr 0 & 0 \cr}\right],
$$
where $Q_1$ is an $n\times n$ idempotent matrix with rank $n-1$.
\end{lemma}

\begin{lemma}\label{kernelimageinclusion}
Let $A, B \in M_{m \times n} (\D)$ be adjacent matrices such that ${\rm rank}\, A = {\rm rank}\, B$. Then
${\rm Im}\, A = {\rm Im}\, B$ or ${\rm Ker}\, A = {\rm Ker}\, B$.
\end{lemma}

\begin{lemma}\label{coadjacency}
Let $A$ and $B$ be matrices of $M_{m \times n}(\D)$, and assume that the rank $r$ of $A$ is greater than $1$.
Suppose in addition that every rank $r-1$ matrix that is adjacent to $A$ is also adjacent to $B$.
Then, either $A=B$, or $B=0$ and $r=2$.
\end{lemma}

{\sl Proof.}
Obviously, no generality is lost in assuming that $A=\left[\matrix{
I_r & 0 \cr
0 & 0}\right]$. First, we prove that $B=\left[\matrix{
C & 0 \cr
0 & 0}\right]$ for some $C \in M_r(\D)$. As at least one rank $r-1$ matrix is adjacent to $A$ (e.g. $\left[\matrix{
I_{r-1} & 0 \cr
0 & 0}\right]$ has this property), we obtain $\rk \,B \in \{r-2,r-1,r\}$. \\
Assume that $\rk\, B<r$.
Let $P \in M_r(\D)$ be an arbitrary rank $r-1$ idempotent. Then, as $Q=\left[\matrix{ P & 0 \cr 0 & 0 }\right]$
has rank $r-1$ and is adjacent to $A$, it is also adjacent to $B$, and hence either $\Im\, B \subset \Im\, Q$ or $\Ker\, Q \subset \Ker\, B$.
If $B$ had rank $r-1$, then we could pick three different $P$'s with distinct kernels and images, and we would obtain a contradiction.
Thus, $B$ has rank $r-2$. Then, $\Im\, B \subset \Im\, Q$ for every matrix $Q$ of the above type, and varying $P$
leads to $B=0$ and hence $r=2$.

Now, we assume that $B \neq 0$ or $r \neq 2$. Then, we have $\rk\, B=r$.

Using again rank $r-1$ idempotents of $M_r(\D)$, we obtain that
$\Im \,Q \subset \Im\, B$ and $\Ker\, B \subset \Ker\, Q$ for any matrix $Q$ as above, and varying $Q$
leads to the fact that $B=\left[\matrix{
B' & 0 \cr
0 & 0
}\right]$ for some invertible matrix $B' \in \GL_r(\D)$. Thus, every rank $r-1$ idempotent $P$ of $M_r(\D)$ is adjacent to $B'$.

We shall finally prove that $B'=I_r$.
We demonstrate first that $B'$ and $I_r$ have the same first $r-1$ rows.
Let us write $B'=\left[\matrix{
B_1 & \, ^{t}y \cr
x & a
}\right]$ with $B_1 \in M_{r-1}(\D)$, $x \in \D^{r-1}$, $^{t}y \in {}^t \D^{r-1}$, and $a \in \D$.
For all $x' \in \D^{r-1}$, the matrix
$\left[\matrix{
I_{r-1} & 0 \cr
x' & 0
}\right]$ is a rank $r-1$ idempotent and hence $\left[ \matrix{
B_1-I_{r-1} & \, ^{t}y \cr
x-x' & a
}\right]$ should have rank $1$.
If $a=0$, then choosing $x' \neq x$ would lead to $^{t}y=0$, which is absurd because the last column of the invertible matrix $B'$ is nonzero.
Thus, $a \neq 0$; then, choosing $x'=x$ leads to $B_1-I_{r-1}=0$, and choosing $x' \neq x$ entails that $^{t}y=0$.

Therefore, $B'$ has the same first $r-1$ rows as $I_r$.
With a similar line of reasoning, one shows that $B'$ has the same last $r-1$ rows as $I_r$, and one concludes that $B'=I_r$.
Therefore, $B=A$ as claimed.
\enp

\subsection{The geometry of adjacent sets for $2 \times 2$ matrices}

Given respective subspaces $V$ and $W$ of $\D^n$ and ${}^t \D^m$, the sets
$V^\bot=\{\, ^{t}x \in {}^t \D^n : \; \forall y \in V, \; y\, ^{t}x=0\}$ and
$W^\bot=\{u \in \D^m : \; \forall \, ^{t}v \in W, \; u \, ^{t}v=0\}$ are subspaces, respectively, of ${}^t \D^n$ and $\D^m$,
with respective dimensions $n-\dim V$ and $m-\dim W$, and we have the relations
$(V^\bot)^\bot=V$ and $(W^\bot)^\bot=W$.

Now, we shall establish a few basic results on the geometry of the maximal adjacent sets in $M_{m \times n}(\D)$.
Most of them deal with the special case $m=n=2$.

The proof of the first lemma is obvious.

\begin{lemma}\label{interlemma}
Let $d$ and $d'$ be points in $\PP(\D^n)$. Then, $R_d \cap R_{d'}=\{0\}$ if and only if $d \neq d'$.  \\
Let $e$ and $e'$ be points in $\PP({}^t\D^m)$. Then, $L_e \cap L_{e'}=\{0\}$ if and only if $e \neq e'$.
\end{lemma}

\begin{lemma}\label{ortholemma}
Let $d \in \PP(\D^2)$ and $e \in \PP({}^t\D^2)$.
Then, $R_d \cap L_e$ contains exactly one rank $1$ idempotent if $d \neq e^\bot$, and it contains
none if $d=e^\bot$.
\end{lemma}

{\sl Proof.}
By definition, $R_d$ is the set of all matrices whose image is included in $d$,
and $L_e$ is the set of all matrices whose kernel includes $e^\bot$.
If $R_d \cap L_e$ contains a rank $1$ idempotent $p$, then the image of $p$ is $d$ and its kernel is $e^\bot$,
and hence $d \neq e^\bot$ and $p$ is the projection of $\D^2$ on $d$ alongside $e^\bot$.
Conversely, if $d \neq e^\bot$ then $\D^2=d \oplus e^\bot$ and the projection $p$ on $d$ alongside $e^\bot$
is a rank $1$ idempotent that belongs to $R_d \cap L_e$.
\enp

\begin{corollary}\label{adjrank2lemma}
Let $A \in M_2(\D)$ be a rank $2$ matrix, and let $d \in \PP (\D ^2)$ and $e \in \PP ({}^t\D^2)$.
Then, $R_d \cap L_e$ contains a rank $1$ matrix that is adjacent to $A$ if and only if
$d \neq (A^{-1} e)^\bot$.
\end{corollary}

{\sl Proof.}
Let $M \in M_2(\D)$ be a rank $1$ matrix. Then, $A^{-1}M$ has rank $1$, and $M$ is adjacent to $A$ if and only if $A^{-1}M$ is adjacent to $I_2$.
On the other hand $M$ belongs to $R_d \cap L_e$ if and only if $A^{-1}M$ belongs to $R_d \cap L_{A^{-1}e}$.
It follows that $R_d \cap L_e$ contains a rank $1$ matrix that is adjacent to $A$ if and only if
$R_d \cap L_{A^{-1}e}$ contains a rank $1$ matrix that is adjacent to $I_2$ which, according to Lemma \ref{subadjacentidempotent},
amounts to the existence of a rank $1$ idempotent in $R_d \cap L_{A^{-1}e}$. The conclusion readily follows from Lemma \ref{ortholemma}.
\enp

Let $d \in \PP(\D^2)$ and $e \in \PP({}^t \D^2)$.
Choosing a nonzero vector $d_1 \in d$, we have a bijection
$${}^ty \in {}^t \D^2 \; \mapsto \; {}^tyd_1 \in R_d,$$
and we can use this bijection to endow $R_d$ with a structure of right vector space over $\D$.
This structure depends on the choice of $d_1$ for if we let $a \in \D \setminus \{0\}$
and denote by $\bullet$ and $\star$ the scalar multiplications associated respectively with $d_1$ and $ad_1$, then
for all $b \in \D$ and ${}^ty \in {}^t \D^2$ we get
$$({}^tyd_1)\star b=({}^ty a^{-1}(ad_1)) \star b=({}^tya^{-1}b)ad_1={}^ty(a^{-1}ba)d_1=({}^tyd_1) \bullet (a^{-1}ba).$$
Nevertheless, this very equality shows that the notion of linear and affine subspaces of $R_d$
does not depend on the choice of the nonzero vector $d_1$ used to define the vector space structure, and neither
does the notion of semilinear maps from or to $R_d$.

In a similar fashion, we can endow $L_e$ with a structure of left vector space over $\D$
by choosing a nonzero vector ${}^te_1 \in e$ and by using the bijection
$$x \in \D^2 \; \mapsto {}^te_1 x \in L_e.$$
The linear and affine subspaces of the resulting vector space are independent from the choice
of the specific vector ${}^te_1$, and so does the notion of semilinear maps from or to $L_e$.

In the following lemma, we characterize the affine lines in $R_d$ and $L_e$ in terms of adjacency:

\begin{lemma}\label{affinelinelemma}
Let $d \in \PP(\D^2)$ and $e \in \PP({}^t \D^2)$.
Let $A$ be a rank $2$ matrix of $M_2(\D)$. Then, the set of matrices of $R_d$ (respectively, of $L_e$)
that are adjacent to $A$ is an affine line that does not go through $0$.
Conversely, given an affine line $\calD$ of $R_d$ (respectively, of $L_e$) that does not go through $0$,
there exists a rank $2$ matrix $B \in M_2(\D)$ such that $\calD$ is the set of all matrices of $R_d$ (respectively, of $L_e$)
that are adjacent to $B$.
\end{lemma}

{\sl Proof.}
Obviously, no generality is lost in assuming that $d=\D \times \{0\}$.
Then, a matrix $M \in R_d$ is adjacent to $A$ if and only if $A^{-1}M$, which belongs to $R_d$, is adjacent to $I_2$, i.e. $A^{-1}M$
is a rank $1$ idempotent of $R_d$. Thus, the set of all matrices of $R_d$ that are adjacent to $A$ is the set of all matrices of the form
$$A \,\left[\matrix{
1 & 0 \cr
a & 0
}\right] \quad \textrm{with $a \in \D$,}$$
and it is obviously an affine line of $R_d$ that does not go through zero.

Conversely, given an affine line $\calD$ of $R_d$ that does not go through zero, we can find
$B \in \GL_2(\D)$ such that $\calD$ is the set of all matrices $B \,
\left[\matrix{
1 & 0 \cr
a & 0
}\right]$ with $a \in \D$, and the above proof shows that $\calD$ is the set of all matrices of $R_d$ that are adjacent to $B$.

The proof is similar in the case of $L_e$.
\enp

\subsection{General results on the action of an adjacency preserver on adjacent sets}

In the next few lemmas we will always assume that $m,n,p,q$ are positive integers and
$\phi : M_{m \times n} (\D) \to M_{p \times q} (\D)$ is an adjacency preserver satisfying $\phi (0) = 0$.
Note that $\phi$ maps rank one matrices to rank one matrices.
The first lemma is well-known. The proof can be found in \cite[ Lemmas 5.3 and 5.4]{Se6}.

\begin{lemma}\label{asqsa}
Let $d$ and $e$ be elements of $\PP(\D^{n})$ and $\PP(^{t} \D^{m})$, respectively.
Then, either there is a unique $d' \in \PP(\D^q)$ such that
$\phi(R_d) \subset R_{d'}$, or there is a unique $e' \in \PP({}^t\D^p)$ such that
$\phi(R_d) \subset L_{e'}$. \\
Moreover, either there is a unique $d' \in \PP(\D^q)$ such that
$\phi(L_e) \subset R_{d'}$, or there is a unique $e' \in \PP({}^t\D^p)$ such that
$\phi(L_e) \subset L_{e'}$.
\end{lemma}

\begin{lemma}\label{pasmolezacel}
Let $d$ and $e$ be elements of $\PP(\D^{n})$ and $\PP(^{t} \D^{m})$, respectively.
If we have some $d' \in \PP(\D^{q})$ such that $\phi(R_d) \subset R_{d'}$ and $\phi(L_e) \subset R_{d'}$,
or we have some $e' \in \PP(^{t} \D^{p})$ such that $\phi(R_d) \subset L_{e'}$ and $\phi(L_e) \subset L_{e'}$,
then $\phi$ is degenerate with respect to $0$.
\end{lemma}

{\sl Proof.}
Assume that we have some $d' \in \PP(\D^{q})$ such that $\phi(R_d) \subset R_{d'}$ and $\phi(L_e) \subset R_{d'}$. \\
Let $d_1 \in \PP(\D^{n})$. If $\phi(R_{d_1}) \subset R_{d_2}$ for some $d_2 \in \PP(\D^q)$, then
as $R_{d_1} \cap L_e$ contains a rank $1$ matrix $A$, we see that $\phi(A)$ has rank $1$ and belongs to
$R_{d_2} \cap R_{d'}$, which yields $d_2=d'$. Thus, either
$\phi(R_{d_1}) \subset R_{d'}$ or $\phi(R_{d_1}) \subset L_{e'}$ for some $e' \in \PP({}^t \D^p)$.
Similarly, for all $e_1 \in \PP({}^t\D^{m})$, either $\phi(L_{e_1}) \subset R_{d'}$ or $\phi(L_{e_1}) \subset L_{e'}$ for some $e' \in \PP({}^t \D^p)$.

Let now $(d_1,e_1) \in \PP(\D^{n}) \times \PP({}^t\D^m)$. Assume that there is a pair
$(e'_1,e'_2) \in \PP({}^t \D^p)^2$ such that $\phi(R_{d_1}) \subset L_{e'_1}$ and $\phi(L_{e_1}) \subset L_{e'_2}$.
Again, as $R_{d_1} \cap L_{e_1}$ contains a rank $1$ matrix and $\phi$ maps rank $1$ matrices to rank $1$ matrices, we deduce that $e'_1=e'_2$.

Now, we can conclude: if $\phi$ maps every rank $1$ matrix into $R_{d'}$ then we are done;
if not then there are points $e_1 \in  \PP({}^t\D^m)$ and $e' \in  \PP({}^t\D^p)$ such that
$\phi(L_{e_1}) \subset L_{e'}$. Then, for every $d_2 \in \PP(\D^n)$, the above proof shows that
$\phi(R_{d_2}) \subset R_{d'}$ or $\phi(R_{d_2}) \subset L_{e'}$, whence $\phi$
maps $M_{m \times n}^{\leq 1}(\D)$ into $L_{e'} \cup R_{d'}$, which shows that $\phi$ is degenerate with respect to $0$.

The second case is handled in a similar fashion.
\enp

\begin{corollary}\label{mainactioncor}
Assume that $\phi$ is not degenerate with respect to $0$.
Then, one of the following two cases holds:
\begin{itemize}
\item[(i)] There are nonconstant maps $\alpha : \PP(\D^n) \rightarrow \PP(\D^q)$ and $\beta :
\PP({}^t\D^m) \rightarrow \PP({}^t\D^p)$ such that
$\phi(R_d) \subset R_{\alpha(d)}$ for all $d \in \PP(\D^n)$ and
$\phi(L_e) \subset L_{\beta(e)}$ for all $e \in \PP({}^t\D^m)$.
\item[(ii)] There are nonconstant maps $\alpha : \PP(\D^n) \rightarrow \PP({}^t\D^p)$ and $\beta :
\PP({}^t\D^m) \rightarrow \PP(\D^q)$ such that
$\phi(R_d) \subset L_{\alpha(d)}$ for all $d \in \PP(\D^n)$ and
$\phi(L_e) \subset R_{\beta(e)}$ for all $e \in \PP({}^t\D^m)$.
\end{itemize}
\end{corollary}

{\sl Proof.}
Fix $d_0 \in \PP(\D^n)$. Assume first that there is a $d'_0 \in \PP(\D^q)$ such that $\phi(R_{d_0}) \subset R_{d'_0}$.
Let $e \in \PP({}^t\D^m)$. If there is a point $d'' \in \PP(\D^q)$ such that $\phi(L_e) \subset R_{d''}$,
then as $R_{d_0} \cap L_e$ contains a rank $1$ matrix we deduce that $R_{d''} \cap R_{d'_0}$ is nonzero and hence $d'_0=d''$,
which by Lemma  \ref{pasmolezacel} would entail that $\phi$ is degenerate with respect to $0$.
Thus, there is a $\beta(e) \in \PP({}^t \D^p)$ such that
$\phi(L_e) \subset L_{\beta(e)}$. With the same line of reasoning, one deduces that for every $d \in \PP(\D^n)$,
there is an $\alpha(d) \in \PP(\D^q)$ such that $\phi(R_d) \subset R_{\alpha(d)}$.
As $\phi$ is not degenerate with respect to $0$, we see that $\alpha$ is nonconstant, and similarly $\beta$
is nonconstant. Thus, situation (i) holds.

If there is an $e'_0 \in \PP({}^t\D^p)$ such that $\phi(R_{d_0}) \subset L_{e'_0}$, then a similar line of reasoning
as the above one shows that situation (ii) holds.
\enp

\subsection{The initial reduction}

\begin{proposition}\label{initialreduction}
Let $m,n,p,q$ be positive integers with $m \geq 2$ and $n\geq 2$, and let
$\phi : M_{m\times n}(\D) \rightarrow M_{p\times q}(\D)$ be an adjacency preserver that maps $0$ to $0$
and is not degenerate with respect to $0$.
Then, there are invertible matrices $P \in \GL_{m}(\D)$ and $Q \in \GL_n(\D)$ such that
the map
$$A \in M_2(\D) \mapsto \phi\left(P \left[\matrix{
A & 0 \cr
0 & 0 }\right] Q\right)$$
is not degenerate with respect to $0$ and maps at least one rank $2$ matrix to a rank $2$ matrix.
\end{proposition}

{\sl Proof.}
As we know that $\phi$ is not degenerate with respect to $0$, it satisfies one of the two possible conclusions of Corollary \ref{mainactioncor}.
We shall assume that condition (i) holds, the proof being similar in the second case.
Then, we have nonconstant maps $\alpha : \PP(\D^n) \rightarrow \PP(\D^q)$ and $\beta :
\PP({}^t\D^m) \rightarrow \PP({}^t\D^p)$ such that
$\phi(R_d) \subset R_{\alpha(d)}$ for all $d \in \PP(\D^n)$ and
$\phi(L_e) \subset L_{\beta(e)}$ for all $e \in \PP({}^t\D^m)$, and hence we can choose
distinct points $d_1,d_2$ of $\PP(\D^n)$ and distinct points $e_1,e_2$ of $\PP({}^t\D^m)$ such that
$\alpha(d_1) \neq \alpha(d_2)$ and $\beta(d_1) \neq \beta(d_2)$.
Obviously, no generality is lost in assuming that $d_1$ and $d_2$ are spanned, respectively, by the first two vectors of the standard basis of $\D^n$,
and ditto for $e_1$ and $e_2$ with respect to the standard basis of ${}^t \D^m$.
Then, we wish to show that the adjacency preserver
$$\theta : A \in M_2(\D) \mapsto \phi\left(\left[ \matrix{
A & 0 \cr
0 & 0
}\right] \right)$$ is not degenerate and maps at least one rank $2$ matrix to a rank $2$ matrix.

Assume then that $\theta$ maps every rank $1$ matrix into $R_d \cup L_e$ for some $d \in \PP(\D^q)$ and some $e \in \PP({}^t \D^p)$.
In particular $\theta(E_{1,1})$ belongs to $R_d \cup L_e$ and hence $d=\alpha(d_1)$ or $e=\beta(e_1)$.
If $d=\alpha(d_1)$, then as $\theta(E_{1,2})$ belongs to $R_{\alpha(d_2)} \cap L_{\beta(e_1)}$
and $\alpha(d_1) \neq \alpha(d_2)$, we obtain $\beta(e_1)=e$. In any case $\beta(e_1)=e$, and in the same manner we obtain $\alpha(d_1)=d$.
Working similarly with $E_{2,2}$ instead of $E_{1,1}$, we would obtain $\alpha(d_2)=d$, contradicting $\alpha(d_1) \neq \alpha(d_2)$.
Thus, $\theta$ is not degenerate with respect to $0$.

Let us now prove that $\theta$ maps at least one rank $2$ matrix to a rank $2$ matrix.
Assume on the contrary that $\theta$ maps every rank $2$ matrix to a matrix of rank at most $1$.
We shall distinguish between two cases.

Assume first that $\card\,\D>2$.
Let $A \in M_2(\D)$ be a rank $2$ matrix with nonzero entries everywhere, so that
$d_i \neq ( A^{-1}e_j)^\bot$ for all $i,j$ in $\{1,2\}$ (indeed, we have $A (d_i)^{\bot} \neq e_j$ for all $(i,j) \in \{1,2\}^2$).
Thus, by Corollary \ref{adjrank2lemma}, for all $i,j$ in $\{1,2\}$, the set $R_{d_i} \cap L_{e_j}$ contains a rank $1$ matrix that is adjacent to $A$,
and hence $R_{\alpha(d_i)} \cap L_{\beta(e_j)}$ contains a rank $1$ matrix that is adjacent to $\theta(A)$.
If $\theta(A)$ had rank $1$, then  $R_{\alpha(d_1)}$ would contain two matrices with distinct kernels
that are adjacent to $\theta(A)$, leading to $\theta(A) \in R_{\alpha(d_1)}$, and similarly one would obtain
$\theta(A) \in R_{\alpha(d_2)}$, a contradiction.
Therefore, $\theta(A)=0$. To conclude, we simply choose $t \in \D \setminus \{0,1\}$
and we note that the matrices $A_1=\left[\matrix{
1 & 1 \cr
t & 1
}\right]$ and $A_2=\left[\matrix{
t & 1 \cr
1 & 1
}\right]$ have rank $2$ and nonzero entries everywhere, leading to $\theta(A_1)=0=\theta(A_2)$,
which contradicts the obvious fact that $A_1$ and $A_2$ are adjacent.

Assume finally that $\D$ has cardinality $2$.
Let $A \in M_2(\D)$ be a rank $2$ matrix. Assume that $\theta(A)$ has rank $1$.
We know from Corollary \ref{adjrank2lemma} that $R_{d_1}$ contains distinct rank $1$ matrices $B_1$ and $B_2$ that are adjacent to $A$;
then, $\theta(B_1)$ and $\theta(B_2)$ are distinct nonzero elements of $R_{\alpha(d_1)}$, which, as $\D$ only has two elements, yields that
they have distinct kernels and their range is $\alpha(d_1)$. As $\theta(A)$ is a rank $1$ matrix that is adjacent to
both of them, we would deduce that $\theta(A) \in R_{\alpha(d_1)}$. Similarly, we obtain that $\theta(A) \in R_{\alpha(d_2)}$,
a contradiction. Thus, we find that $\theta$ maps every rank $2$ matrix to $0$, and once more we obtain a contradiction by considering
two adjacent rank $2$ matrices (this time, we can take $\left[\matrix{
1 & 0 \cr
0 & 1
}\right]$ and $\left[\matrix{
1 & 1 \cr
0 & 1
}\right]$).

It follows that at least one rank $2$ matrix is mapped by $\theta$ to a rank $2$ matrix,
which completes the proof.
\enp

\subsection{The key lemma}

In this section, we shall prove the most important point of our proof:

\begin{proposition}\label{keylemma}
Let $p,q \geq 2$ be integers, $\D$ be an EAS division ring, and
$\phi : M_2(\D) \rightarrow M_{p\times q}(\D)$ be an adjacency preserver such that $\phi(0)=0$.
Assume also that $\phi$ does not map the set of rank $1$ matrix into an adjacent set and that $\phi$ maps at least one rank $2$ matrix to a rank $2$ matrix.
Then, $\phi$ is standard.
\end{proposition}

The proof has many steps. We will successively compose $\phi$ on the right and on the left with several standard adjacency preservers
that map $0$ to $0$ so as to ultimately find a map of the form
$M \mapsto \left[\matrix{
M & 0 \cr
0 & 0
}\right]$. Note that composing $\phi$ with such maps alters none of our basic assumptions.

\vskip 3mm
\noindent \textbf{Step 1.} An initial reduction. \\
We choose a rank $2$ matrix $M_0 \in M_2(\D)$ such that $\phi(M_0)$ has rank $2$.
Thus, we can write $\phi(M_0)=P \left[\matrix{
I_2 & 0 \cr
0 & 0
}\right]Q$ for invertible matrices $P \in \GL_p(\D)$ and $Q \in \GL_q(\D)$.
Replacing $\phi$ with $M \mapsto P^{-1}\phi(M_0 M) Q^{-1}$, we are effectively reduced to the case when
$\phi(I_2)=\left[\matrix{
I_2 & 0 \cr
0 & 0}\right]$.

\vskip 3mm
\noindent \textbf{Step 2.} The action of $\phi$ on rank $1$ idempotent matrices. \\
Let $M \in M_2(\D)$ be a rank $1$ idempotent. Thus $M$ has rank $1$ and is adjacent to $I_2$,
whence $\phi(M)$ has rank $1$
and is adjacent to $\phi(I_2)$. We deduce from Lemma \ref{subadjacentidempotent} that
$$\phi(M)=\left[\matrix{
K(M) & 0 \cr
0 & 0
}\right] \quad \textrm{for some rank $1$ idempotent $K(M)\in M_2(\D)$.}$$
In particular, we can find an invertible matrix $R \in \GL_2(\D)$ such that
$$\phi(E_{1,1})=\left[\matrix{
RE_{1,1} R^{-1} & 0 \cr
0 & 0
}\right].$$
Setting $P:=R \oplus I_{p-2}$ and $Q:=R \oplus I_{q-2}$, and replacing $\phi$ with
$M \mapsto P^{-1}\phi(M)Q^{-1}$, we see that no further generality is lost in assuming that
$$\phi(E_{1,1})=E_{1,1}.$$
Note here that the $E_{1,1}$ on the left-hand side of the equality belongs to $M_2(\D)$, whereas the one on the right-hand side belongs to $M_{p\times q}(\D)$.

\vskip 3mm
\noindent \textbf{Step 3.} $\phi$ is not degenerate with respect to $0$. \\
Assume on the contrary that there are points $d \in \PP(\D^q)$ and $e \in \PP({}^t \D^p)$
such that $\phi$ maps every rank $1$ matrix into $R_d \cup L_e$.

Suppose that there exist three distinct points $d_1,d_2,d_3$ in $\PP(\D^2)$ such that
$\phi(R_{d_i}) \subset R_d$ for all $i=1,2,3$.
As $\phi$ does not map the set of rank $1$ matrices into $R_d$, we can find $e_1 \in \PP({}^t\D^2)$
such that $\phi(L_{e_1}) \subset L_e$.
Without loss of generality, we may then assume that $d_1 \neq (e_1)^\bot$ and $d_2 \neq (e_1)^\bot$,
and hence we can find rank $1$ idempotent matrices $A_1 \in R_{d_1} \cap L_{e_1}$ and $A_2 \in R_{d_2} \cap L_{e_1}$.
As $\phi(A_1)$ and $\phi(A_2)$ must then belong to $R_d \cap L_e$, the matrices
$K(A_1)$ and $K(A_2)$ are rank $1$ idempotents with the same kernel and range, and hence they are equal,
contradicting the fact that $A_1$ and $A_2$ are adjacent.

With the same line of reasoning, we obtain that there cannot exist three distinct points
$d'_1,d'_2,d'_3$ in $\PP(\D^2)$ such that $\phi(R_{d'_i}) \subset L_e$ for all $i=1,2,3$.
If $\D$ contains at least $4$ elements, then this is a contradiction because $\PP(\D^2)$ contains at least $5$ points.

Assume finally that $\D$ is finite with cardinality denoted by $k$. Then, the above arguments show that we can choose
distinct points $d_1$ and $d_2$ in $\PP(\D^2)$ such that
$\phi(R_{d_i}) \subset R_d$ for $i=1,2$, or $\phi(R_{d_i}) \subset L_e$ for $i=1,2$.
Let us assume that the first case holds (the second one is handled similarly). Again,
we can find a point $e_1 \in \PP({}^t\D^2)$ such that $\phi(L_{e_1}) \subset L_e$.
Then, $\phi$ is injective on $L_{e_1}$ and maps every element of
$L_{e_1} \cap (R_{d_1} \cup R_{d_2})$ into $R_d \cap L_e$.
This is absurd because the first set has cardinality $2k-1$ and the second one has cardinality $k$.

We conclude that $\phi$ is not degenerate with respect to $0$.
Thus, we can apply Corollary \ref{mainactioncor}. In what follows, we will only consider the case
when condition (ii) is satisfied, that is when we have nonconstant maps
$\alpha : \PP(\D^2) \rightarrow \PP({}^t\D^p)$ and $\beta :
\PP({}^t\D^2) \rightarrow \PP(\D^q)$ such that
$\phi(R_d) \subset L_{\alpha(d)}$ for all $d \in \PP(\D^2)$ and
$\phi(L_e) \subset R_{\beta(e)}$ for all $e \in \PP({}^t\D^2)$.
The case when condition (i) holds can be handled is a similar manner.

\vskip 3mm
\noindent \textbf{Step 4.} The maps $\alpha$ and $\beta$ are injective. \\
Assume that there are distinct points $d$ and $d'$ in $\PP(\D^2)$
such that $\alpha(d)=\alpha(d')$. The projective line $\PP({}^t \D^2)$ has at least three points so we can choose
one such point $e$ that is different from $d^\bot$ and $(d')^\bot$, so that
we can find rank $1$ idempotents $A \in R_d \cap L_e$ and $A' \in R_{d'} \cap L_e$.
Then, $A$ and $A'$ are adjacent as they are distinct elements of $L_e$, and hence
$K(A)$ and $K(A')$ are distinct rank $1$ idempotents. However, $K(A)$ and $K(A')$ would have the same range and the same
kernel since $\phi(A)$ and $\phi(A')$ are both included in $L_{\alpha(d)} \cap R_{\beta(e)}$. This is absurd.

We conclude that $\alpha$ is injective, and the same line of reasoning shows that $\beta$ is also injective.
In particular, we deduce that $\phi$ maps rank $1$ matrices with distinct images to rank $1$ matrices with distinct kernels,
and maps rank $1$ matrices with distinct kernels to rank $1$ matrices with distinct images.

\vskip 3mm
\noindent \textbf{Step 5.} On the ranges of $\alpha$ and $\beta$. \\
Let $d \in \PP(\D^2)$. By Lemma \ref{ortholemma}, we can choose a rank $1$ idempotent $A \in R_d$, and hence
$K(A)$ is a rank $1$ idempotent while $\phi(A) \in L_{\alpha(d)}$. It follows that $\alpha(d) \in \PP({}^t (\D^2 \times \{0\}))$.
Similarly, one proves that the range of $\beta$ is included in $\PP(\D^2 \times \{0\})$.
In particular, for every rank $1$ matrix $A$, we have
$$\phi(A)=\left[\matrix{
K(A) & 0 \cr
0 & 0
}\right] \quad \textrm{for some $K(A) \in M_2(\D)$.}$$

\vskip 3mm
\noindent \textbf{Step 6.} The action of $\phi$ on rank $2$ matrices. \\
In this step, we prove that $\phi$ maps every rank $2$ matrix to a rank $2$ matrix. Let
$A$ be a rank $2$ matrix in $M_2(\D)$. Assume firstly that $\phi(A)$ has rank $1$.
Let $d \in \PP(\D^2)$. We know from Corollary \ref{adjrank2lemma} that $R_d$ contains two (rank $1$) matrices $B_1$ and $B_2$ that are adjacent to
$A$ and have distinct kernels. Then, $\phi(B_1)$ and $\phi(B_2)$ are rank $1$ matrices of $L_{\alpha(d)}$ that are adjacent to $\phi(A)$
and have distinct ranges. Therefore, $\phi(A)$ must belong to $L_{\alpha(d)}$.
Varying $d$ yields a contradiction as $\alpha$ is injective.

Thus, $\phi$ must map each rank $2$ matrix to a matrix of rank $0$ or $2$. Assume then that some rank $2$ matrix $A\in M_2(\D)$
satisfies $\phi(A)=0$. Obviously, the rank $2$ matrix $A':=A \, \left[\matrix{
1 & 1 \cr
0 & 1
}\right]$ is adjacent to $A$ and hence $\phi(A')$ should have rank $1$, contradicting the above result.

Therefore, $\phi$ maps rank $2$ matrices to rank $2$ matrices.

\vskip 3mm
\noindent \textbf{Step 7.} Reduction to the case when $p=q=2$. \\
Let $A \in M_2(\D)$ be of rank $2$. Let $d \in \PP(\D^2)$.
We can choose a rank $1$ matrix $B$ in $R_d$ that is adjacent to $A$, and hence
$\phi(B)$ is a rank $1$ matrix that is adjacent to the rank $2$ matrix $\phi(A)$, leading to
$\Ker \,\phi(A) \subset \Ker\, \phi(B)=\alpha(d)^\bot$. As $\alpha$ is injective and has its range included in $\PP({}^t (\D^2 \times \{0\}))$, varying $d$
yields $\Ker\, \phi(A) \subset \{0\} \times  \D^{p-2}$, whence
$\Ker\, \phi(A) = \{0\} \times  \D^{p-2}$.
Similarly, one shows that $\Im \,\phi(A)=\D^2 \times \{0\}$.

Thus, we have shown that, for all $M \in M_2(\D)$, we have some $K(M) \in M_2(\D)$
with the same rank as $M$ and such that
$$\phi(M)=\left[\matrix{
K(M) & 0 \cr
0 & 0
}\right].$$
Replacing $\phi$ with $K$, we are reduced to the following situation:
$p=q=2$, $\phi$ fixes $I_2$ and $E_{1,1}$, and $\phi$ preserves adjacency, idempotency and the rank.

\vskip 3mm
\noindent \textbf{Step 8.} The action of $\phi$ on the $R_d$ and $L_e$ spaces. \\
Let $d \in \PP(\D^2)$ and $e \in \PP({}^t \D^2)$.
Let us prove that $\phi$ induces a lineation from $R_d$ to $L_{\alpha(d)}$.
We already know that $\phi$ takes $1$-dimensional linear subspaces of $R_d$ to $1$-dimensional linear subspaces of $L_{\alpha(d)}$.
Indeed, collinear nonzero vectors of $R_d$ have the same kernel and hence $\phi$ takes them to matrices of $L_{\alpha(d)}$
with the same image, and hence to collinear vectors of $L_{\alpha(d)}$.
Next, let $\calD$ be an affine line of $R_d$ that does not go through zero. Then, Lemma \ref{affinelinelemma} yields
a rank $2$ matrix $A\in M_2(\D)$ such that $\calD$ is the set of all matrices of $R_d$ that are adjacent to $A$.
Therefore, $\phi(\calD)$ is included in the set of all matrices of $L_{\alpha(d)}$ that are adjacent to the rank $2$ matrix $\phi(A)$,
an affine line by Lemma \ref{affinelinelemma}. Thus, $\phi$ induces a lineation from $R_d$ to $L_{\alpha(d)}$, as claimed.
Its image is not contained in an affine line since it contains $0$ and $\phi$ maps matrices with different kernels
to matrices with different images. Finally, we already knew that $\phi$ is injective on $R_d$ and that $\phi(0)=0$.

Thus, if $\card \,\D>2$ then by the fundamental theorem of affine geometry $\phi$ induces a semilinear isomorphism
from $R_d$ to $L_{\alpha(d)}$; if $\card\, \D=2$ the result is also true simply because it is known in this case that a map
between $2$-dimensional spaces over $\D$ is linear provided that it is bijective and maps $0$ to $0$.

In the same way, we obtain that $\phi$ induces a semilinear isomorphism from $L_e$ to $R_{\beta(e)}$.

\vskip 3mm
\noindent \textbf{Step 9.} Reduction to a map that induces the identity on $R_{d_1}$. \\
From now on we set
$$d_1:=\D \times \{0\} \subset \D^2.$$
Note that $\alpha(d_1)=e_1:={}^td_1$ since $\phi(E_{1,1})=E_{1,1}$. It follows from Step 8 that the
restriction of $\phi$ to $R_{d_1}$ can be written as
$M \mapsto {}^t(M^\tau)Q$ for some invertible matrix $Q \in \GL_2(\D)$ and some antiautomorphism $\tau$ of $\D$.
As $\phi(E_{1,1})=E_{1,1}$, we have
$Q=\left[\matrix{
1 & 0 \cr
a & b
}\right]$ for some $(a,b)\in \D^2$. Moreover, as the rank $1$ idempotent matrix $P=\left[\matrix{
1 & 0 \cr
1 & 0
}\right]$ of $R_{d_1}$ is mapped by $\phi$ to a rank $1$ idempotent, and on the other hand
${}^t (P^\tau) Q=\left[\matrix{
1 + a & b \cr
0 & 0
}\right]$, we must have $1+a=1$ and hence $a=0$.
It follows that $N \mapsto Q^{-1}N$ is the identity on $L_{e_1}$,
whence $\phi$ coincides with $M \mapsto Q^{-1}\, {}^t(M^\tau) Q$ on $R_{d_1}$.

Now, we set
$$\psi : M \in M_2(\D) \mapsto {}^t(Q\phi(M)Q^{-1})^{\tau^{-1}}.$$
With the above results, we obtain the following properties of $\psi$:
\begin{itemize}
\item $\psi$ preserves adjacency, idempotency, and the rank;
\item $\psi$ is the identity on $R_{d_1}$;
\item There are injections $\gamma : \PP(\D^2) \rightarrow \PP(\D^2)$ and $\eta : \PP({}^t\D^2) \rightarrow \PP({}^t\D^2)$
such that $\psi$ induces semilinear bijections from $R_d$ to $R_{\gamma(d)}$ and from $L_e$ to $L_{\eta(e)}$
for every $d$ in $\PP(\D^2)$ and every $e$ in $\PP({}^t \D^2)$.
\end{itemize}

In the next steps, we shall prove that $\psi$ is the identity on the whole $M_2(\D)$, and it will follow that $\phi$ is standard, as claimed.

\vskip 3mm
\noindent \textbf{Step 10.} The maps $\gamma$ and $\eta$ are identities.

Let $e \in \PP({}^t\D^2)$. We see that $\psi$ maps $R_{d_1} \cap L_e$ (which contains a rank $1$ matrix) into $R_{d_1} \cap L_{\eta(e)}$.
However, $\psi$ fixes every point in $R_{d_1}$ whence $L_e$ has a nonzero common matrix with $L_{\eta(e)}$. Therefore, $\eta(e)=e$.
Thus, $\eta$ is the identity on $\PP({}^t \D^2)$.

We have just proved that
$\psi$ induces a semilinear automorphism of $L_e$. Then, we deduce that
$\gamma$ is surjective: indeed, given $d' \in \PP(\D^2)$, we know that $L_e \cap R_{d'}$ contains a nonzero point and we can then choose
a nonzero $A \in L_e$ such that $\psi(A) \in L_e \cap R_{d'}$. However, $A \in R_{d}$ for some $d \in \PP(\D^2)$
and hence $R_{d'}$ and $R_{\gamma(d)}$ have a nonzero common matrix, which proves that $\gamma(d)=d'$.
Therefore, $\gamma$ is a bijection.

Let $d \in \PP(\D^2)$.
Then, for all $d' \in \PP(\D^2) \setminus \{d\}$, we know that
$R_{d'} \cap L_{d^\bot}$ contains a rank $1$ idempotent, whence its image by $\psi$ is a rank $1$ idempotent in $R_{\gamma(d')} \cap L_{d^\bot}$
, which yields $\gamma(d') \neq d$.
Thus, $\gamma$ is a permutation of $\PP(\D^2)$ that maps $\PP(\D^2) \setminus \{d\}$ into itself, and the conclusion follows that $\gamma(d)=d$.

It follows that $\psi$ preserves both the kernel and the image of rank $1$ matrices.

\vskip 3mm
\noindent \textbf{Step 11.} $\psi$ leaves every rank $1$ matrix invariant.

Let $P \in M_2(\D)$ be a rank $1$ idempotent. Then, $\psi(P)$ is a rank $1$ idempotent with the same kernel and image
as $P$, and hence $\psi(P)=P$.

Let $d \in \PP(\D^2)$. We know that $\psi$ induces a semilinear automorphism of $R_d$.
However, by Lemma \ref{affinelinelemma} the set of rank $1$ idempotents in $R_d$ is an affine line that does not go through zero,
and $\psi$ leaves all its points invariant. By looking at $\psi$ on the translation vector space of this line,
we see that $\psi$ is actually linear (i.e., the automorphism of $\D$
that is associated to $\psi$ is the identity), and then as the above affine line spans the vector space $R_d$ we conclude that
$\psi$ is the identity on $R_d$. Varying $d$ yields that $\psi$ fixes every matrix of $M_2^{\leq 1}(\D)$.

\vskip 3mm
\noindent \textbf{Step 12.} $\psi$ is the identity.

We already know that $\psi$ fixes every matrix with rank at most $1$.
Let $A$ be a rank $2$ matrix. Then, every rank $1$ matrix that is adjacent to $A$
is also adjacent to the rank $2$ matrix $\psi(A)$. Using Lemma \ref{coadjacency}, we conclude that $\psi(A)=A$.
Therefore, $\psi$ is the identity and we conclude that $\phi$ is standard.

\subsection{Completing the reduction to the standard case}

Throughout this section, the elements of the standard basis of the left vector space
$\D^n$ (respectively, the right vector space
$^{t}\D^m$) will be denoted by $e_1 , \ldots , e_n$
(respectively, by $^{t}f_1 , \ldots, \, ^{t}f_m$).

Our aim here is to prove the following general result:

\begin{proposition}\label{extensionproposition}
Let $\D$ be an EAS division ring,
$m,n,p,q \ge 2$ be integers, and $\phi : M_{m\times n}(\D) \rightarrow M_{p \times q}(\D)$ be an adjacency preserver.
Assume that
$$\phi(A)=\left[\matrix{
A & 0 \cr
0 & 0
}\right] \quad \textrm{for every $A \in M_2(\D)$.}$$
Then, $m \leq p$, $n \leq q$, and $\phi$ is standard.
\end{proposition}

The proof is split into three lemmas, the first of which is a consequence of
Proposition \ref{keylemma}.

\begin{lemma}\label{partialstandard}
Let $\D$ be an EAS division ring,
$p,q \ge 2$ be integers, and $\phi : M_2 (\D ) \to M_{p \times q} (\D)$ be an adjacency preserver.
Assume that there exists a nonzero scalar $\alpha \in \D$ such that
\begin{equation}\label{mkurc}
\phi \left( \left[ \matrix{ a & 0 \cr b & 0 \cr}\right] \right) =  \left[ \matrix{  \left[ \matrix{ a & 0 \cr b & 0 \cr}\right] & 0 \cr 0 & 0 \cr}\right]
\quad \textrm{for every $a,b \in \D$,}
\end{equation}
\begin{equation}\label{mkurc2}
\phi \left( \left[ \matrix{ 0 & 1 \cr 0 & 0 \cr}\right] \right) =  \left[ \matrix{  \left[ \matrix{ 0 & 1 \cr 0 & 0 \cr}\right] & 0 \cr 0 & 0 \cr}\right] \ \ \ {\rm and}\ \ \
\phi \left( \left[ \matrix{ 0 & 0 \cr 0 & 1 \cr}\right] \right) =  \left[ \matrix{  \left[ \matrix{ 0 & 0 \cr 0 & \alpha \cr}\right] & 0 \cr 0 & 0 \cr}\right].
\end{equation}
Then
$$
\phi (A) =  \left[ \matrix{  A & 0 \cr 0 & 0 \cr}\right] \quad \textrm{for every $A\in M_2 (\D)$.}
$$
\end{lemma}

{\sl Proof.} We see that $\phi(0)=0$.
It is obvious from identities (\ref{mkurc}) and (\ref{mkurc2}) that $\phi$ is not degenerate with respect to $0$.
Thus, Propositions \ref{initialreduction} and \ref{keylemma} entail that $\phi$ is standard, that is
either (\ref{onopr}), or (\ref{onodr}) with $R=0$. In the second case we would have $\phi ( R(e_1)) \subset L({}^te)$ for some
${}^te \in {}^t \D^p$, in contradiction with
(\ref{mkurc}). Hence, $\phi$ is of the form
 (\ref{onopr}), and since $\phi ( E_{1,1}) = E_{1,1}$ and $\phi (E_{2,2}) = \alpha E_{2,2}$ we have
$$
\phi (A) =  \left[ \matrix{  \left[ \matrix{ t_1 & 0 \cr 0 & t_2  \cr}\right]\, A^\tau \,  \left[\matrix{ s_1 & 0 \cr 0 & s_2 \cr}\right] & 0 \cr 0 & 0 \cr}\right] \quad \textrm{for every $A \in M_2 (\D)$.}
$$
Here, $t_1 s_1= 1$ and $t_2s_2 = \alpha$.

In particular, $\phi$ maps every matrix $A\in M_2 (\D)$ into a $p \times q$ matrix having nonzero entries only in the upper left $2\times 2$ corner. Of course, we can forget the bordering zeroes, or in other words, we may consider $\phi$ as a map of $M_2 (\D)$ into itself.

From $\phi ( a E_{1,1}) = aE_{1,1}$ we get that $t_1 \tau (a) t_1^{-1} = a$ for every $a\in \D$, whence
$\tau : a \mapsto t_1^{-1} a t_1$. It follows that $A^\tau=(t_1^{-1}I_2)A(t_1I_2)$ for every $A \in M_2(\D)$.
Setting $t= t_2 t_1^{-1}$ and $s= t_1 s_2$, we deduce that
$$
\phi (A) =    \left[ \matrix{ 1 & 0 \cr 0 & t \cr}\right]\, A \,  \left[\matrix{ 1 & 0 \cr 0 & s \cr}\right] \quad \textrm{for every $A \in M_2(\D)$.}
$$
The equations $\phi (E_{1,2}) = E_{1,2}$ and $\phi (E_{2,1}) = E_{2,1}$ yield the desired conclusion that $s=t=1$.
\enp

\begin{lemma}\label{jlikoli}
Let $m,n,p,q \ge 2$ be integers, $\D$ be an EAS division ring, and
$\phi : M_{m \times n} (\D) \to M_{p \times q} (\D)$ be an adjacency preserver.
Assume that
$$
\phi \left(   \left[ \matrix{ A & 0 \cr 0 & 0 \cr }\right] \right) =   \left[ \matrix{ A & 0 \cr 0 & 0 \cr }\right]
\quad \textrm{for every $A \in M_2 (\D)$.}
$$
Then, $p \geq m$, $q \geq n$ and there are invertible matrices $T \in \GL_p(\D)$ and $S \in \GL_q(\D)$
such that $\phi$ coincides with $M \mapsto T \left[\matrix{ M & 0 \cr 0 & 0 }\right] S$
on every matrix of $R(e_1) \cup L({}^tf_1)$.
\end{lemma}

\bigskip

{\sl Proof.}
We will use an induction process.
First, given $r \in \{ 2,\dots,n-1\}$, we prove that if there  exists an invertible matrix
$$
V= \left[ \matrix{ I_2 & 0 \cr V_1 & V_2 \cr} \right] \in M_q (\D)
$$
such that $q \geq r$ and
$$
\phi \left(   \left[ \matrix{ A & 0 \cr 0 & 0 \cr }\right] \right) =   \left[ \matrix{ A & 0 \cr 0 & 0 \cr }\right] \, V \quad \textrm{for every rank $1$ matrix $A \in M_{2\times r} (\D)$,}
$$
then $q \ge r+1$ and there exists an invertible matrix
$$
S= \left[ \matrix{ I_2 & 0 \cr S_1 & S_2 \cr} \right] \in M_q (\D)
$$
such that
$$
\phi \left(   \left[ \matrix{ A & 0 \cr 0 & 0 \cr }\right] \right) =   \left[ \matrix{ A &  0 \cr  0 & 0 \cr }\right] S \quad \textrm{for  every rank $1$ matrix $A \in M_{2\times (r+1)} (\D)$.}
$$
When proving this induction step we can assume with no loss of generality that $V=I_q$.

Thus, assume that $q \geq r$ and
\begin{equation}\label{marak}
\phi \left(   \left[ \matrix{ A & 0 \cr 0 & 0 \cr }\right] \right) =   \left[ \matrix{ A & 0 \cr 0 & 0 \cr }\right]
\quad \textrm{for every rank $1$ matrix $A \in M_{2\times r} (\D)$.}
\end{equation}
Then clearly $\phi (L (\, ^{t}f_1)) \subset L(\, ^{t}f_1)$. It follows that
$\phi (E_{1, r+1}) = \sum_{j=1}^q \lambda_j E_{1,j}$ for some scalars $\lambda_1 , \ldots , \lambda_q$.
If $\lambda_{r+1} = \ldots = \lambda_q = 0$, then we would have
$$
\phi \left(  \sum_{j=1}^r \lambda_j E_{1,j} \right) =  \sum_{j=1}^r \lambda_j E_{1,j} = \phi (E_{1, r+1}),
$$
contradicting the fact that $ \sum_{j=1}^r \lambda_j E_{1,j}$ and $E_{1, r+1}$ are adjacent.
Hence, at least one of scalars $\lambda_j$, $j \ge r+1$, is nonzero. In particular, $q \ge r+1$.
It follows that there exists an invertible  matrix $S \in M_q (\D)$ such that
$$
e_j S = e_j, \ \ \ j = 1, \ldots , r, \ \ \ {\rm and} \ \ \  \left(  \sum_{j=1}^q \lambda_j e_j \right) S = e_{r+1}.
$$
Clearly, every such $S$ has the block matrix form
$$
S = \left[ \matrix{ I_r & 0 \cr S_1 & S_2 \cr}\right].
$$
Replacing $\phi$ by the map $A\mapsto \phi (A) S$, $A\in M_{m \times n} (\D)$, we may, and we will assume that
we have (\ref{marak}), and additionally, $\phi (E_{1, r+1}) = E_{1, r+1}$.
We now use the fact that the rank one matrix $\phi ( E_{2, r+1}) \in \phi (L (\, ^{t}f_2)) \subset L( \, ^{t}f_2 )$
is adjacent to  $\phi (E_{1, r+1}) = E_{1, r+1}$ to conclude that there exists a nonzero scalar $\alpha$ such
that $\phi ( E_{2, r+1}) = \alpha E_{2, r+1}$.

Our next goal is to show that if $u \in {\rm span}\, \{ e_1 , \ldots , e_r \} \subset \D^n$ is a nonzero vector, then
for all scalars $\alpha_{11}, \alpha_{12}, \alpha_{21}, \alpha_{22} \in \D$ we have
$$
\phi (\, ^{t} f_1 \alpha_{11} u + \, ^{t} f_1 \alpha_{12} e_{r+1} +  \, ^{t} f_2 \alpha_{21} u  + \, ^{t} f_2 \alpha_{22} e_{r+1} )
$$
\begin{equation}\label{karzetako}
= \, ^{t} f_1 \alpha_{11} u + \, ^{t} f_1 \alpha_{12} e_{r+1} +  \, ^{t} f_2 \alpha_{21} u  + \, ^{t} f_2 \alpha_{22} e_{r+1}
\end{equation}
(note that $^{t}f_1$ and $^{t}f_2$ in the first line of the equation denote vectors in $\, ^{t}\D^m$, while the same symbols in
the second line stand for the first two canonical basic vectors in $\, ^{t}\D^p$, and that similarly, the vectors $u$ and $e_{r+1}$ in the first line
are not the same vectors as the vectors $u$ and $e_{r+1}$ appearing in the second line).
All we need to do is to observe that the set of all matrices $\, ^{t} f_1 \alpha_{11} u + \, ^{t} f_1 \alpha_{12} e_{r+1} +  \, ^{t} f_2 \alpha_{21} u  + \, ^{t} f_2 \alpha_{22} e_{r+1}$,
$\alpha_{11}, \alpha_{12}, \alpha_{21}, \alpha_{22} \in \D$, can be identified with $M_2 (\D)$.
Then (\ref{karzetako}) follows directly from Lemma \ref{partialstandard}.

It follows that we have
$$
\phi \left(   \left[ \matrix{ A & 0 \cr 0 & 0 \cr }\right] \right) =   \left[ \matrix{ A & 0 \cr 0 & 0 \cr }\right]
$$
for every $A \in M_{2 \times (r+1)} (\D)$ of rank one.

Thus, by induction we find that $n \leq q$ and that there is an invertible matrix
$$
S= \left[ \matrix{ I_2 & 0 \cr S_1 & S_2 \cr} \right] \in \GL_q (\D)
$$
such that
$$
\phi \left(   \left[ \matrix{ A \cr 0}\right] \right) =   \left[ \matrix{ A & 0 \cr 0 & 0 \cr }\right]\,S \quad \textrm{for every rank $1$ matrix $A \in M_{2\times n} (\D)$.}
$$
With the same line of reasoning, we find that $m \leq p$ and that there is an
invertible matrix
$$
T= \left[ \matrix{ I_2 & T_1 \cr 0 & T_2 \cr} \right] \in \GL_p (\D)
$$
such that
$$
\phi \left(   \left[ \matrix{ A & 0 }\right] \right) =   T\,\left[ \matrix{ A & 0 \cr 0 & 0 \cr }\right] \quad \textrm{for every rank $1$ matrix $A \in M_{m\times 2} (\D)$.}
$$
With the respective shapes of $S$ and $T$, it is then obvious that
$\phi$ coincides with $M \mapsto T \left[\matrix{ M & 0 \cr 0 & 0 }\right] S$
on every matrix of $R(e_1) \cup L({}^tf_1)$.
\enp

\begin{lemma}\label{jlikoli3}
Let $m,n,p,q \ge 2$ be integers with $m \leq p$ and $n \leq q$, let $\D$ be an EAS division ring, and
$\phi : M_{m \times n} (\D) \to M_{p \times q} (\D)$ be an adjacency preserver that is not degenerate with respect to $0$.
For $M \in M_{m \times n}(\D)$, set
$$\widetilde{M}:=\left[\matrix{
M & 0 \cr
0 & 0
}\right] \in M_{p \times q}(\D).$$
Assume that $\phi(M)=\widetilde{M}$ for all $M \in R(e_1) \cup L({}^t f_1)$. \\
Then, $\phi(M)=\widetilde{M}$ for all $M \in M_{m \times n}(\D)$.
\end{lemma}

{\sl Proof.}
We shall prove that $\phi$ coincides with $\psi : M \mapsto \widetilde{M}$
first on rank $1$ matrices, and then on all matrices by induction on the rank.

It is obvious from our assumptions that $\phi$ cannot satisfy conclusion (ii) from Corollary \ref{mainactioncor}, and hence it satisfies conclusion (i),
yielding maps $\alpha : \PP(\D^n) \rightarrow \PP(\D^q)$ and $\beta : \PP({}^t\D^m) \rightarrow \PP({}^t\D^p)$
such that $\phi(R_d) \subset R_{\alpha(d)}$ and $\phi(L_e) \subset L_{\beta(e)}$ for all $d \in \PP(\D^n)$ and all $e \in \PP({}^t\D^m)$.
We naturally identify $\PP(\D^n)$ with a subspace of $\PP(\D^q)$, and ditto for
$\PP({}^t\D^m)$ with respect to $\PP({}^t\D^p)$.
Given $d \in \PP(\D^n)$, we can choose a nonzero $A \in R_d \cap L({}^tf_1)$, and hence
$\phi(A)=A$, which yields $\alpha(d)=d$. Similarly, one finds that $\beta$ is the identity on $\PP({}^t\D^m)$.

Now, let $B$ be a rank $1$ matrix of $M_{m \times n} (\D)$ that does not belong to $R(e_1) \cup L({}^tf_1)$.
Without loss of generality, we can assume that $B \in R(e_2) \cap L({}^tf_2)$.
Then, as $\alpha$ and $\beta$ are identities, we find $\phi(B) \in R(e_2) \cap L({}^tf_2)$, and of course $\phi(B) \neq 0$.
Then, we see that the map
$A \in M_2(\D) \mapsto \phi\left(\left[\matrix{
A & 0 \cr
0 & 0
}\right]\right)$ satisfies all the assumptions of Lemma \ref{partialstandard}, and we conclude in particular that $\phi(B)=\widetilde{B}$.
Thus, we have shown that $\phi$ coincides with $\psi$ on all rank $1$ matrices.

From there, we proceed by induction on the rank.
Let $k \geq 1$, and assume that $\phi$ coincides with $\psi$ on all matrices of rank at most $k$.
Let $B \in M_{m\times n}(\D)$ be a rank $k+1$ matrix. Let $M$ be a rank $k$ matrix that is adjacent to $\widetilde{B}$.
Then, $\Im M \subset \Im \widetilde{B}$ and $\Ker \widetilde{B} \subset \Ker M$, and hence
$M=\widetilde{M'}$
for some matrix $M' \in M_{m\times n}(\D)$, and obviously $M'$ is adjacent to $B$.
By induction, we know that $M=\phi(M')$, and hence $M$ is adjacent to $\phi(B)$.
Using Lemma \ref{coadjacency}, we conclude that $\phi(B)=\psi(B)$, or $\phi(B)=0$ and $k=1$.
Now, if $k=1$ and $\phi(B)=0$ we deduce that $\phi(B')$ has rank $1$ for every rank $2$ matrix that is adjacent to $B$,
and this is absurd as we have just shown that $\phi$ maps no rank $2$ matrix to a rank $1$ matrix
(and at least one rank $2$ matrix is adjacent to $B$).

Thus, by induction we conclude that $\phi=\psi$, as claimed.
\enp

Now, we complete the proof of Proposition \ref{extensionproposition}.
Let $\phi$ be as in Proposition \ref{extensionproposition}. Obviously $\phi$
is not degenerate with respect to $0$.
Then, Lemma \ref{jlikoli} yields that $m \leq p$, $n \leq q$ and that there are invertible matrices $S \in \GL_q(\D)$
and $T \in \GL_p(\D)$ such that $\psi : M \mapsto T \phi(M)S$ satisfies the assumptions of Lemma
\ref{jlikoli3}; then, one concludes that $\phi(M)=T^{-1} \left[\matrix{
M & 0 \cr
0 & 0
}\right]S^{-1}$ for all $M \in M_{m \times n}(\D)$, and hence $\phi$ is standard.

\vskip 3mm
\noindent {\em Remark:} In Proposition \ref{extensionproposition}, the assumption that the division ring $\D$ be EAS
is unnecessary. To prove this, one simply needs to generalize Lemma \ref{partialstandard}
to all division rings. We shall not do this as it is not directly useful to our study.

\subsection{Wrapping the proof up}

Now, we can conclude. First, we obtain the following result.

\begin{proposition}\label{conclusionprop}
Let $m,n,p,q$ be positive integers, and $\phi : M_{m\times n}(\D) \rightarrow M_{p\times q}(\D)$ be an adjacency preserver such that $\phi(0)=0$
and $\phi$ is not degenerate with respect to $0$. Then, $\phi$ is standard.
\end{proposition}

{\sl Proof.}
The integers $m,n,p,q$ are all greater than $1$ since otherwise $\phi$ would be degenerate with respect to $0$.
By Proposition \ref{initialreduction}, we lose no generality in assuming that there is a map $\theta : M_2(\D) \rightarrow M_{p \times q}(\D)$
that preserves adjacency, maps $0$ to $0$, is not degenerate with respect to $0$, maps at least one rank $2$ matrix to a rank $2$ matrix, and satisfies
$$\phi \left( \left[ \matrix{ A & 0 \cr 0 & 0}\right]\right)=\theta(A), \quad \textrm{for every $A \in M_2(\D)$.}$$
Then, Proposition \ref{keylemma} yields that $\theta$ is standard. Without further loss of generality, we may then assume that
$$\theta(A)=\left[\matrix{
A & 0 \cr
0 & 0
}\right], \quad \textrm{for every $A \in M_2(\D)$.}$$
Thus, Proposition \ref{extensionproposition} applies to $\phi$ and shows that it is standard.
\enp

From there, we can easily complete the proof of Theorem \ref{maintheorem}:
let $\phi : M_{m\times n}(\D) \rightarrow M_{p\times q}(\D)$ be a non-standard adjacency preserver.
Let $A \in M_{m\times n}(\D)$. If $\phi$ were not degenerate with respect to $A$,
then $\psi : N \mapsto \phi(A+N)-\phi(A)$ would be an adjacency preserving map that maps $0$ to $0$
and is not degenerate with respect to $0$, and hence Proposition \ref{conclusionprop} would show that $\psi$ is standard.
Then, as $\phi(M)=\phi(A)+\psi(M-A)$ for all $M \in M_{m\times n}(\D)$, the conclusion would follow that
$\phi$ is standard as well, contradicting our assumptions.
Thus, $\phi$ is degenerate.

\section{Degenerate adjacency preservers over finite fields}

Our goal in this section is to prove Theorem \ref{finitefieldstheo}.
Throughout the section, $\F$ is an arbitrary finite field whose cardinality is denoted by $c$.

\subsection{Preliminary lemmas}

\begin{lemma}\label{cardinalitylemma1}
Let $A$ be a rank $r$ matrix of $M_{m \times n}(\F)$, with $r \geq 2$. Then,
there is an adjacent set $\mathcal{V}$ with cardinality $c^{r-1}$ in which all the
matrices have rank $r-1$ and are adjacent to $A$.
\end{lemma}

{\sl Proof.}
By left and right-multiplying $A$ with invertible matrices, we see that no generality is lost in assuming
that $A=\left[\matrix{
I_r & 0 \cr
0 & 0
}\right]$. Then, one simply takes the set $\mathcal{V}$ of all matrices of the form
$\left[\matrix{
I_{r-1} & {}^tx & 0 \cr
0 & 0 & 0
}\right]$ with ${}^tx \in {}^t \F^{r-1}$.
\enp

\begin{lemma}\label{cardinalitylemma2}
Let $A$ be a rank $2$ matrix of $M_{m \times n}(\F)$.
Let $d \in \PP(\F^n)$ and $e \in \PP({}^t\F^m)$.
Then, the set of matrices in $R_d$ that are adjacent to $A$ is empty or has exactly $c$ elements, and only
the second option is possible if $m=n=2$; the same is true of the set of matrices in $L_e$ that are adjacent to $A$. \\
Moreover, the set of rank $1$ matrices in $A+R_d$ is empty or has exactly $c$ elements, and only the second option is possible if $m=n=2$;
the same holds for the set of rank $1$ matrices in $A+L_e$.
\end{lemma}

{\sl Proof.}
As in the proof of the previous lemma, no generality is lost in assuming that $A=\left[\matrix{
I_2 & 0 \cr
0 & 0
}\right]$.

If $d \not\subset \F^2 \times \{0\}$ - which cannot occur if $n=2$ - then we deduce from the first principles given
in the beginning of Section \ref{basicresults}
that no matrix of $R_d$ is adjacent to $A$. If $d \subset \F^2 \times \{0\}$ then
Lemma \ref{subadjacentidempotent} yields a bijection between the set of matrices of $R_d$
that are adjacent to $A$ and the set of rank $1$ idempotents that belong to $R_{d'} \subset M_2(\F)$, where
$d'$ is the point of $\PP(\F^2)$ that corresponds to $d$ through $x \in \F^2 \mapsto (x,0) \in \F^q$.
Thus, by Lemma \ref{affinelinelemma}, $R_d$ contains exactly $c$ matrices that are adjacent to $A$.

On the other hand $M \mapsto A-M$ maps bijectively the set of all rank $1$ matrices of $A+R_d$
onto the set of all matrices of $R_d$ that are adjacent to $A$, whence the first one has at most $c$ elements.

One proceeds in a similar way to prove the other two statements of the lemma.
\enp

\subsection{The key result}

\begin{proposition}\label{keylemmafinite}
Let $p,q$ be integers with $p \geq 2$ and $q \geq 2$, let $\F$ be a finite field, and
$\phi : M_2(\F) \rightarrow M_{p \times q}(\F)$ be a degenerate
adjacency preserver such that $\phi(0)=0$. Then,
$\phi(A)$ has rank at most $1$ for all $A \in M_2(\F)$.
\end{proposition}

{\sl Proof.}
The proof has several steps.
Assume that $\phi$ maps some rank $2$ matrix $A$ to a rank $2$ matrix.

First of all, since $\phi$ is degenerate, it is not standard and we deduce from Proposition \ref{keylemma}
that its maps the set of rank $1$ matrices of $M_2(\F)$ into an adjacent set of rank $1$ matrices. Then, one of the following two situations holds:
\begin{itemize}
\item \textbf{Case (a):} there is a point $d \in \PP(\F^q)$ such that $\phi$ maps every rank $1$ matrix into $R_d$;
 \item \textbf{Case (b):} there is a point $e \in \PP({}^t\F^p)$ such that $\phi$ maps every rank $1$ matrix into $L_e$.
\end{itemize}
We shall only tackle Case (a) as the proof is essentially similar in Case (b).
\vskip 3mm
\noindent \textbf{Step 1.} $\phi$ maps every rank $2$ matrix that is adjacent to $A$ to a rank $2$ matrix. \\
First of all, we note that $M \mapsto \phi(A)-\phi(A-M)$ satisfies the same assumptions as above
and hence our initial step yields that either there is a point $d' \in \PP(\F^q)$ such that $\phi({\mathcal B} (A,1)) \subset
\phi(A)+R_{d'}$ or there is a point $e' \in \PP({}^t\F^p)$ such that
$\phi({\mathcal B} (A,1)) \subset \phi(A)+L_{e'}$. In any case, Lemma \ref{cardinalitylemma2}
shows that $\phi({\mathcal B} (A,1))$ contains at most $c$ rank $1$ matrices.
Now, let $d_1 \in \PP(\F^2)$. Then, Lemma \ref{cardinalitylemma2} yields that
$A+R_{d_1}$ contains $c$ rank $1$ matrices, and as $\phi$ is injective on it we deduce that
 $\phi$ maps every rank $2$ matrix of $A+R_{d_1}$
to a rank $2$ matrix. Varying $d_1$ yields that $\phi$ maps every rank $2$ matrix of ${\mathcal B} (A,1)$ to a rank $2$ matrix.
\vskip 3mm
\noindent \textbf{Step 2.} $\phi$ preserves the rank. \\
Let $B \in M_2(\F)$ be of rank $2$. Let $a \in \F^*$ and $b \in \F \setminus \{0,1\}$. Whenever
$P$ equals either one of the elementary operation matrices $\left[\matrix{1 & 0 \cr 0 & b}\right]$,
$\left[\matrix{1 & a \cr 0 & 1}\right]$ or $\left[\matrix{1 & 0 \cr a & 1}\right]$,
the matrix $PB$ has rank $2$ and is adjacent to $B$.
By Gaussian elimination we know that, for every rank $2$ matrix $B$,
we can find matrices $P_1,\dots,P_N$ of the above type such that $B=P_1\cdots P_N A$, and by induction
we deduce from Step 1 that $\phi(B)$ has rank $2$. Thus, we have shown that $\phi$ preserves the rank.

\vskip 3mm
\noindent \textbf{Step 3.} $\phi$ is injective and preserves adjacency in both directions. \\
Let $B$ be a rank $1$ matrix of $M_2(\F)$. We choose a rank $2$ matrix $A_1 \in M_2(\F)$
that is not adjacent to $B$. Let $d_1 \in \PP(\F^2)$ be such that $B \in R_{d_1}$.
Then, $R_{d_1}$ contains $c$ matrices that are adjacent to $A_1$, and $\phi$ maps $R_{d_1}$ injectively
into $R_d$, which contains at most $c$ matrices that are adjacent to $\phi(A_1)$. It follows
that $\phi(B)$ is not adjacent to $\phi(A_1)$.

Now, let $B \in M_2(\F)$. We can find a matrix $B' \in M_2(\F)$
such that $\phi(B)$ and $\phi(B')$ are non-adjacent: if $B$ has rank $1$ this is deduced from the above proof;
if $B$ has rank $2$, we take $B'=0$; if $B=0$ we take $B'=A$.
Then, $\psi : M \mapsto \phi(M+B)-\phi(B)$ is a degenerate adjacency preserver
that maps $0$ to $0$ and the rank $2$ matrix $B'-B$ to a rank $2$ matrix. Applying Step 2 to $\psi$ shows that
for every $B'' \in M_2(\F) \setminus \{B\}$ that is non-adjacent to $B$, the matrix
$\phi(B'')$ is different from $\phi(B)$ and not adjacent to it.

Therefore, $\phi$ is injective and preserves adjacency in both directions.

\vskip 3mm
\noindent \textbf{Step 4.} The final contradiction. \\
The matrices $E_{1,1}$ and $E_{2,2}$ are non-adjacent and have rank $1$.
By Step 3, the matrices $\phi(E_{1,1})$ and $\phi(E_{2,2})$ should be distinct and non-adjacent.
On the other hand, they must belong to $R_d$. This is absurd.

Thus, our initial assumption was wrong and $\phi$ must map every rank $2$ matrix to a matrix with rank at
most $1$. This completes the proof.
\enp

\subsection{Completing the proof}

\begin{proposition}\label{extensionlemmafinite}
Let $m,n,p,q$ be positive integers, $\F$ be a finite field,
and $\phi : M_{m \times n}(\F) \rightarrow M_{p \times q}(\F)$ be a degenerate adjacency preserver
such that $\phi(0)=0$. Then, $\phi$ maps every matrix to a matrix with rank at most $1$.
\end{proposition}

{\sl Proof.}
The result is obvious if one the integers $m,n,p,q$ equals $1$, so in the rest of the proof we assume that this is not the case.

The proof works by induction on the rank.
We already know that $\phi$ maps every rank $1$ matrix to a rank $1$ matrix.
Next, let $A$ be a rank $2$ matrix. Without loss of generality we can assume that $A=\left[\matrix{I_2 & 0 \cr
0 & 0}\right]$. Then,
$$\theta : M \in M_2(\F) \mapsto \phi\left(\left[
\matrix{M & 0 \cr 0 & 0}\right]\right)$$
is a degenerate adjacency preserver that maps $0$ to $0$. By Proposition \ref{keylemmafinite},
$\theta(I_2)$ has rank at most $1$, i.e.\ $\phi(A)$ has rank at most $1$.

From there, we proceed by induction. Let $r \geq 3$ be such that $r \leq m$ and $r \leq n$,
and assume that $\phi$ maps every rank $r-1$ matrix to a matrix with rank at most $1$.
Let $A$ be a rank $r$ matrix. Choosing a rank $r-1$ matrix $B$ that is adjacent to $A$,
we deduce that $\rk\, \phi(A) \leq \rk\, \phi(B)+1 \leq 2$. Assume then that $\rk  \, \phi(A)=2$.
Lemma \ref{cardinalitylemma1} helps us find
an adjacent set $\mathcal{V}$ in $M_{m \times n}(\F)$ with cardinality $c^{r-1}$ in which all the matrices have rank $r-1$ and are adjacent to $A$.
Then, $\phi(\mathcal{V})$ must be an adjacent set
in $M_{p \times q}(\F)$ in which all the matrices have rank at most $1$, and hence
it is included in $R_d$ for some $d \in \PP(\F^q)$ or in $L_e$ for some
$e \in \PP({}^t \F^p)$; as all the matrices of $\phi(\mathcal{V})$ must be adjacent to $\phi(A)$,
we deduce from Lemma \ref{cardinalitylemma2} that $\phi(\mathcal{V})$ has cardinality at most $c$,
contradicting the fact that $\phi$ must be injective on $\mathcal{V}$. This is absurd, and we deduce that $\rk \, \phi(A) \leq 1$.

Thus, by induction we conclude that $\phi$ maps $M_{m\times n}(\F)$ into $M_{p \times q}^{\leq 1}(\F)$.
\enp

Now, we can complete the proof of Theorem \ref{finitefieldstheo}.
Let $\phi : M_{m\times n}(\F) \rightarrow M_{p\times q}(\F)$ be a degenerate adjacency preserver.
Let $A \in M_{m \times n}(\F)$. The map $\psi : N \in M_{m \times n}(\F) \mapsto \phi(A+N)-\phi(A)$
satisfies the conditions of Proposition \ref{extensionlemmafinite} and hence $\rk\,(\phi(M)-\phi(A)) \leq 1$ for all $M \in M_{m \times n}(\F)$. Thus, $\phi\left(M_{m \times n}(\F)\right)$ is an adjacent set.

\end{document}